%%%%%%%%%%%%%%  Algebraic and Geometric Topology: agt-2-24.tex  %%%%%%%%
%%%%        
%%%%      A note on the Lawrence-Krammer-Bigelow representation
%%%%             
%%%%                  Luisa Paoluzzi and Luis Paris  
%%%%                          
%%%%                Published in Volume 2(2002) 499--518
%%%%
%%%%                    Publication date 25 June 2002
%%%%
%%%%                      This is a plain TeX file
%%%%
%%%%
%%%%%%%%%%%%%%%%%%                                   %%%%%%%%%%%%%%%%%%%

%%%%%%%%%%%%%%%%%%%%%%%%%%%%%%%%%%%%%%%%%%%%%%%%%%%%%%%%%%%%%%%
%%%%%%%%%%%             gtmacros.tex            %%%%%%%%%%%%%%%
%%%%%%%%%%%             version 1.6             %%%%%%%%%%%%%%% 
%
%                       Colin Rourke   
%
%
%    These macros are recommended for use by authors submitting articles   
%    to Geometry and Topology or to Algebraic and Geometric Topology.  
%    They are intended to be used with plain TeX. Each macro is described 
%    briefly to make it clear how to use it (or to modify it to achieve
%    different results).  If you modify this file then please change its
%    name.  If you modify this file and use the modified file to 
%    format an article for submission to Geometry and Topology or
%    Algebraic and Geometric Topology, then please paste the modified
%    file into your main TeX file.  Do not submit it as a separate file.
%      
%    Instructions on using these macros are also given in  gtmacins.tex  
%    or  gtmacins.ps  or .pdf  available on the gt www pages or by 
%    anonymous ftp from the gt/info/macros directory.
%
%
\magnification=\magstephalf      % Sets default point size to 11pt.
%
%  Basic layout parameters :
%
\vsize=7.5truein                 % Sets text height to 7.5 inches.
\hsize=5.2truein                 % Sets text width to 5.2 inches.
\newskip\stdskip                 % standard vertical space
\stdskip=6pt plus3pt minus3pt    % (slightly more stretchy
\medskipamount=\stdskip          % than the usual \medskip)
\parindent=0pt                   % Paragraphs are non-indented with
\parskip=\stdskip                % a little space between paragraphs. 
\abovedisplayskip=\stdskip       %  Reduces the space
\belowdisplayskip=\stdskip       %  around displays.
\mathsurround=0.75pt             % Gives a little extra space around maths.
\overfullrule=0pt                %  Prevents black boxes
%
%   The following macro is for principal paragraph breaks ie
%   a paragraph break with a slightly larger space :
%
\def\ppar{\par\goodbreak\vskip 8pt plus 4pt minus 4pt}     
%
%  The standard horizontal space for theorems, labels etc :
%
\def\stdspace{\hskip 0.75em plus 0.15em\ignorespaces}
\let\qua\stdspace % useful abbreviation (3/4 of a quad)
%
%%%%%%%%%%%%%%            FONT MACROS            %%%%%%%%%%%%
%
%           The following font macros define the AMS symbol 
%           and Euler-Fraktal fonts for use in text and
%           mathematics with appropriate size changes.
%           They also define two new control sequences  
%           \small  and  \large  (similar to those built
%           into LaTeX) which change the size of all fonts 
%           both in text and maths.  \small  is 10% smaller 
%           than normal and  \large  30% bigger.  The strange
%           size of the \large text fonts (10pt scaled 1315)
%           is because these macros are intended to be used
%           at \magstephalf.  The result is 10pt scaled 1440
%           (\magstep2) which is a standard font size.  If
%           you are borrowing these macros to use them at
%           another basic  \magnification, then you will
%           probably need to change 1315 to 1200 in the eleven
%           places marked ** below.  \large  will then be
%           20% bigger than normal.  Note that at \magstephalf
%           all the fonts come out roughly one point larger
%           than their size as defined in these macros.
%
%           The size-changing macros are based on Knuth's
%           \ninepoint and \eightpoint macros.
%
%
%    The macros are laid out in a way which makes it clear how to
%    add futher fonts (or delete unavailable ones) and how to add
%    further size changes.
%
%    First comes a definition of  \hexnumber  which is needed to
%    refer to font families whose family number is not known :
%
\def\hexnumber#1{\ifcase#1 0\or 1\or 2\or 3\or 4\or 5\or 6\or 7\or 8\or
 9\or A\or B\or C\or D\or E\or F\fi}
%
%     Next we define the AMS symbol-a fonts at 13,10,9,7,6,5 points
%
\font\thirtnmsa=msam10 scaled 1315    %%% **  see note above 
\font\tenmsa=msam10          \font\ninemsa=msam9
\font\sevenmsa=msam7         \font\sixmsa=msam6
\font\fivemsa=msam5
%%%%%%  (add further sizes here if you need them)
%
%    and the standard size family for these fonts
%
\newfam\msafam                  \textfont\msafam=\tenmsa
\scriptfont\msafam=\sevenmsa    \scriptscriptfont\msafam=\fivemsa
\edef\hexa{\hexnumber\msafam}        %  The msa family is  \fam\hexa
\def\msa{\fam\msafam\tenmsa}         %  \msa  switches to this family
%
%    Repeat these steps for the AMS symbol-b fonts
%
\font\thirtnmsb=msbm10 scaled 1315   %%%  ** see note above
\font\tenmsb=msbm10      \font\ninemsb=msbm9
\font\sevenmsb=msbm7     \font\sixmsb=msbm6
\font\fivemsb=msbm5
%%%%%%  (add further sizes here if you need them)
%
\newfam\msbfam                   \textfont\msbfam=\tenmsb       
\scriptfont\msbfam=\sevenmsb     \scriptscriptfont\msbfam=\fivemsb
\edef\hexb{\hexnumber\msbfam}    %  The msb family is \fam\hexb  
\def\msb{\fam\msbfam\tenmsb}     %  \msb switches to this family
%
%        Repeat for the Euler-Fraktal fonts 
%
\font\thirtneufm=eufm10 scaled 1315   %%% **  see note above 
\font\teneufm=eufm10                 \font\nineeufm=eufm9
\font\seveneufm=eufm7                \font\sixeufm=eufm6
\font\fiveeufm=eufm5
%%%%%%  (add further sizes here if you need them)
%
\newfam\eufmfam                    \textfont\eufmfam=\teneufm
\scriptfont\eufmfam=\seveneufm     \scriptscriptfont\eufmfam=\fiveeufm
\edef\hexf{\hexnumber\eufmfam}      % The Euler-Fraktal family is
\def\frak{\fam\eufmfam\teneufm}     % \fam\hexf and \frak switches to this
%
%%%  Add further fonts families here (using the same format) if you need
%    them.  The def of hexnumber is optional (it is only used for
%    \mathchardef 's).
%
%      Now we need to define the standard fonts (which are
%      already defined at 10,7 and 5 point) at 13,9 and 6 point:
%
%      Roman fonts:
\font\thirtnrm=cmr10 scaled 1315    %%%  ** see note above
\font\ninerm=cmr9                   \font\sixrm=cmr6   
%%%%%%  (add further sizes here if you need them)
%
%      Math italic fonts
\font\thirtni=cmmi10 scaled 1315    %%%  ** see note above 
\font\ninei=cmmi9                   \font\sixi=cmmi6  
%%%%%%  (add further sizes here if you need them)
%
%     Symbol fonts
\font\thirtnsy=cmsy10 scaled 1315   %%%  ** see note above
\font\ninesy=cmsy9                  \font\sixsy=cmsy6  
%%%%%%  (add further sizes here if you need them)
%
%     Bold face
\font\thirtnbf=cmbx10 scaled 1315   %%%  ** see note above 
\font\ninebf=cmbx9                  \font\sixbf=cmbx6  
%%%%%%  (add further sizes here if you need them)
%
%     The maths extension font (only defined at text size)
%
\font\thirtnex=cmex10 scaled 1315   %%%  ** see note above
\font\nineex=cmex9                  
%%%%%%  (add further sizes here if you need them)
%
%     Finally three fonts (text italic, slanted and typewriter type)
%     which are also only defined at text size
%
\font\thirtnit=cmti10 scaled 1315  %%%  ** see note above 
\font\nineit=cmti9                  
%%%%%%  (add further sizes here if you need them)
%
\font\thirtnsl=cmsl10 scaled 1315  %%%  ** see note above 
\font\ninesl=cmsl9                  
%%%%%%  (add further sizes here if you need them)
%
\font\thirtntt=cmtt10 scaled 1315  %%%  ** see note above 
\font\ninett=cmtt9                  
%%%%%%  (add further sizes here if you need them)
%
%
%     Now come the two main macros.  What  \small  does is to
%     change all the families of fonts from normal size which is
%     10,7,5  (ie 10pt text, 7pt subscript, 5pt subsubscript)
%     to 9,6,5.  \large  similarly changes to  13,9,7.  To make
%     other size changing macros, choose your three sizes, add
%     font size definitions if necessary and make the obvious changes
%     to one of these macros.  Change  \normalbaselineskip  and
%     \strutbox  dimensions to appropriate sizes as well.  To
%     add further fonts, insert them in each macro, using the
%     AMS fonts as a model.
%      
%
\def\small{%
%
%   redefine the sizes of the roman fonts :
%
\textfont0=\ninerm \scriptfont0=\sixrm \scriptscriptfont0=\fiverm
\def\rm{\fam0\ninerm}%       % ( \rm  sets \ninerm  in text mode
%                            %  and \fam0 in math mode)
%
%   and the math italic fonts :
%
\textfont1=\ninei \scriptfont1=\sixi \scriptscriptfont1=\fivei
%
%   and the symbol fonts :
%
\textfont2=\ninesy \scriptfont2=\sixsy \scriptscriptfont2=\fivesy
%
%   There is only one math extension font :
%
\textfont3=\nineex \scriptfont3=\nineex \scriptscriptfont3=\nineex
%
%   Next the bold font (named rather than numbered) :
%
\textfont\bffam=\ninebf \scriptfont\bffam=\sixbf
\scriptscriptfont\bffam=\fivebf \def\bf{\fam\bffam\ninebf}%
%
%   and the three text-only fonts : 
%
\textfont\itfam=\nineit \def\it{\fam\itfam\nineit}%
\textfont\slfam=\ninesl \def\sl{\fam\slfam\ninesl}%
\textfont\ttfam=\ninett \def\tt{\fam\ttfam\ninett}%
%
%   Now the three new families of AMS fonts :
%
%   AMS symbol-a
%
\textfont\msafam=\ninemsa \scriptfont\msafam=\sixmsa
\scriptscriptfont\msafam=\fivemsa \def\msa{\fam\msafam\ninemsa}%         
%
%   AMS symbol-b
%
\textfont\msbfam=\ninemsb \scriptfont\msbfam=\sixmsb
\scriptscriptfont\msbfam=\fivemsb \def\msb{\fam\msbfam\ninemsb}%         
%
%   Euler-Fraktal font
%
\textfont\eufmfam=\nineeufm  \scriptfont\eufmfam=\sixeufm
\scriptscriptfont\eufmfam=\fiveeufm \def\frak{\fam\eufmfam\nineeufm}%
%
%%%  Add further fonts families here if you need them.
%
%    Reset \normalbaselineskip and \strubox :
%
\normalbaselineskip=11pt%
\setbox\strutbox=\hbox{\vrule height8pt depth3pt width0pt}%
%
%    Set \normalbaselines and \rm (roman) as defaults :
%
\normalbaselines\rm
%
%    Reset some of the basic vertical skips:
%
\stdskip=4pt plus2pt minus2pt    
\medskipamount=\stdskip          
\parskip=\stdskip                
\abovedisplayskip=\stdskip       
\belowdisplayskip=\stdskip       
\def\ppar{\par\goodbreak\vskip 6pt plus 3pt minus 3pt}%     
%
%   And finally reset the size of section heads (see below):
%
\def\section##1{\global\advance\sectionnumber by 1
\vskip-\lastskip\penalty-800\vskip 20pt plus10pt minus5pt 
\egroup{\bf\number\sectionnumber\quad##1}\bgroup\small         
\vskip 6pt plus3pt minus3pt
\nobreak\resultnumber=1}%      % Reset resultnumber at start of section
}    %%%   End of  \small  macro      
%
%   Two useful abbreviations to keep track of \small material:
\def\beginsmall{\bgroup\small}
\let\endsmall\egroup
%
%
%    The \large  macro is similar (comments abbreviated):
%
%
\def\large{%
\textfont0=\thirtnrm \scriptfont0=\ninerm \scriptscriptfont0=\sevenrm
\def\rm{\fam0\thirtnrm}%
\textfont1=\thirtni \scriptfont1=\ninei \scriptscriptfont1=\seveni
\textfont2=\thirtnsy \scriptfont2=\ninesy \scriptscriptfont2=\sevensy
\textfont3=\thirtnex \scriptfont3=\thirtnex \scriptscriptfont3=\thirtnex
\textfont\bffam=\thirtnbf \scriptfont\bffam=\ninebf
\scriptscriptfont\bffam=\sevenbf \def\bf{\fam\bffam\thirtnbf}%
\textfont\itfam=\thirtnit \def\it{\fam\itfam\thirtnit}%
\textfont\slfam=\thirtnsl \def\sl{\fam\slfam\thirtnsl}%
\textfont\ttfam=\thirtntt \def\tt{\fam\ttfam\thirtntt}%
%   AMS symbol-a  :
\textfont\msafam=\thirtnmsa \scriptfont\msafam=\ninemsa
\scriptscriptfont\msafam=\sevenmsa \def\msa{\fam\msafam\thirtnmsa}%         
%   AMS symbol-b  :
\textfont\msbfam=\thirtnmsb \scriptfont\msbfam=\ninemsb
\scriptscriptfont\msbfam=\sevenmsb \def\msb{\fam\msbfam\thirtnmsb}%         
%   Euler-Fraktal font :
\textfont\eufmfam=\thirtneufm  \scriptfont\eufmfam=\nineeufm
\scriptscriptfont\eufmfam=\seveneufm \def\frak{\fam\eufmfam\teneufm}%
%%%% Add further fonts families here if you need them.
%   Reset \normalbaselineskip and \strubox and initialise :
\normalbaselineskip=16pt%
\setbox\strutbox=\hbox{\vrule height11.5pt depth4.5pt width0pt}%
\normalbaselines\rm}%
\let\Large\large   %  for compatibility with latex
%
%   The next two lines define commonly used switches for
%   blackboard bold (\Bbb) and gothic type (\goth).  The   
%   \Bbb  switch is set to work in the same way as in amstex
%   and switches only the next character to blackboard bold.
\def\Bbb#1{{\msb#1}}

%
%   To use the new AMS fonts you can either use the control
%   sequences \msa, \msb (alias \Bbb) and \frak (alias \goth) eg :

   % see the msam font table
%
%   or, more generally, make \mathchardef's (cf Knuth p155) eg :
\mathchardef\plussquare="0\hexa01
\mathchardef\nge="3\hexb0B
\mathchardef\maltesecross="0\hexa7A
\mathchardef\del="0\hexf01
%
%   or you can use the amstex names for all the new symbols by
%   inserting the line  \input amsnames  in your file directly
%   after \input gtmacros. 
%   This presupposes that you have collected a copy of the file
%   amsnames.tex  from the  gt/info/macros  ftp directory.
%
%
%   Finally we need a small capital font (for author(s)) :
%
\font\sc=cmcsc10
%
%%%%%%%%%%%%%%%%%       END OF FONT MACROS     %%%%%%%%%%%%%
%
%
%                 Knuth's \square macro :
%
\def\sqr#1#2{{\vcenter{\vbox{\hrule  height.#2truept
	\hbox{\vrule width.#2truept height#1truept 
	\kern#1truept \vrule width.#2truept}
	\hrule height.#2truept}}}}
\def\sq{\sqr55}    %   A small square for end-of-proofs. 
%                  %   (Define other size squares by varing the
%                  %   the two numbers.)
%
%
%      Style macros for section heads, theorem statements etc :
%   
%
\newcount\sectionnumber            %%%  Allocate registers to take
\newcount\resultnumber             %%%  section and result numbers.
\sectionnumber=0\resultnumber=1    %%%  Set these registers to 0 and 1
%
%   The \section macro produces a \large bold faced section heading
%   numbered to the left.  Pagebreaks are encouraged before the
%   start of the section and discouraged directly after the heading.
%   Typical use  \section{First steps}  with typical result :
%
%    1  First Steps     (set bold and \large)
%
\def\section#1{\global\advance\sectionnumber by 1
\xdef\nextkey{\number\sectionnumber}%      (used by the \key macro)
\vskip-\lastskip\penalty-800\vskip 20pt plus10pt minus5pt 
{\large\bf\number\sectionnumber\quad#1}         
\vskip 8pt plus4pt minus4pt
\nobreak\resultnumber=1}      % Reset resultnumber at start of section
%
%
%
%   Next a macro to set subheadings (like the  \section  macro
%   but without the number, with less space and set in standard size).
%
%   Typical use :  \sh{Example formats}
%
         
%
%   The \proc ... \endproc macros ("proclaim") are for setting theorems, 
%   lemmas, conjectures etc with automatic numbering.  Typical use :    
%  
%    \proc{Theorem}Every lemon is yellow.\endproc
%
%   Typical result :
%     
%    Theorem 3.4  Every lemon is yellow.   

%   (with Theorem 3.4 set bold and a \stdspace of space before the 
%   statement set in slanted type).
%
\def\proc#1{\xdef\nextkey{\number\sectionnumber.\number\resultnumber}%
\vskip-\lastskip\ppar\bf%
\noindent#1\ \number\sectionnumber.\number\resultnumber
\stdspace\sl\global\advance\resultnumber by 1\ignorespaces}
\def\endproc{\rm\ppar} 
%
%  The \prf ... \endprf macros are for setting proofs.  The code
%  for \prf includes the code for \endproc, so there is no need to
%  type \endproc if the theorem is followed immediatedly by a proof.
%
\def\prf{\vskip-\lastskip\ppar\noindent{\bf Proof}%
\stdspace\rm}                            %  For start of proofs  
   %  For end (or absence) of proofs
\def\endprf{\unskip\stdspace\hbox{}%     %  For end of proof (with
\hfill$\sq$\par\medskip}                 %  extra vertical space)  
\def\proof#1{\vskip-\lastskip\ppar\noindent{\bf#1}%
\stdspace\rm\ignorespaces}        %  For start of proof with alternative name
              %  \endproof is an alias for \endprf
%
%   Typical uses :    
%  
%    \proc{Theorem}Every lemon is yellow. \qed\endproc
%
%    \proc{Theorem}Every lemon is yellow.
%    \prf Use your eyes. \endprf
%
%    \proc{Theorem}Every lemon is yellow.
%    \proof{Proof of theorem} Use your eyes. \endprf
%
%   The next macro is a variant of the \proc macro.  It has
%   exactly the same result except that it omits the number.
%
%   Typical use :  
%    
%    \proclaim{Conjecture}Some oranges are yellow.\endproc
%

%
%   The next macro is a further variant for remarks, definitions etc.   
%   It omits the number and does not switch on slanted type.  
%  
%   Typical use :
%
%    \rk{Remark}Some lemons are thick-skinned.\endrk
%

%
%   The next macro is for numbering equations etc, \label  produces the 
%   correct number  x.y  and advances the resultnumber register
%
%   Typical use :
%
%     $$fx=7\eqno{\bf\label}$$
%
%   result :
%
%                           fx = 7                           3.5
%
\def\label{\xdef\nextkey{\number\sectionnumber.\number\resultnumber}%
\number\sectionnumber.\number\resultnumber
\global\advance\resultnumber by 1}
%
%
%
%   The next macros are to automate external references.  To use them 
%   type \reflist ..... \endreflist near the beginning of your paper, 
%   where  .... is the list of references in alphabetical order 
%   and in  the form  \key{KEY}  reference    where "KEY" is a 
%   string of characters which reminds you of the reference.   
%   Separate  references with a blank line or a \par.   Eg 
%
%     \reflist
%
%     ..... more references ....
%
%     \key{Kn-84} {\bf D Knuth}, {\it The TeXbook}, Addison--Wesley (1984)
%
%     ..... more references ....
%
%     \endreflist
%
%   Then type  \references  where you wish the references to be printed
%   (normally near the end of the paper).  To refer to Knuth type
%   for example    see Knuth [\ref{Kn-84}, page 320]   and the correct
%   numerical reference will be printed.  Edit the \references macro
%   to change the formatting (if desired).
%   There is an alternative \refkey for \key, provided your KEY contains
%   only letters.  The syntax is:
%
%     \reflist
%
%     ..... more references ....
%
%     \refkey\Knuth  {\bf D Knuth}, {\it The TeXbook}, Addison--Wesley (1984)
%
%     ..... more references ....
%
%     \endreflist
%
%   \key{Knuth}  has exactly the same maening as \refkey\Knuth and you
%   can mix the two syntaxes if you want.  But \refkey\Kn-84
%   would not work.  It would set Kn as the KEY and -84 would get printed!
%
\newcount\refnumber              %  Register for reference numbers
\refnumber=1                     %  set initially to 1.
\long\def\reflist#1\endreflist{%
\long\def\thereflist{#1}{\def\refkey##1##2\par{\xdef##1{\number\refnumber}%
\global\advance\refnumber by 1}%
\def\key##1##2\par{\expandafter\xdef%
\csname##1\endcsname{\number\refnumber}%
\global\advance\refnumber by 1}#1\par}}
\long\def\references{%
\penalty-800\vskip-\lastskip\vskip 15pt plus10pt minus5pt 
{\large\bf References}\ppar %`References' is set \large bold with space around.
{\leftskip=25pt\frenchspacing    % The list of references is set 
\small\parskip=3pt plus2pt       % \small  with small spaces between,
\def\refkey##1##2\par{\noindent  % numbers in [,]'s and set just to the
\llap{[##1]\stdspace}\ignorespaces##2\par}         % left of a 25pt margin.
\def\key##1##2\par{\noindent  
\llap{[\ref{##1}]\stdspace}\ignorespaces##2\par}  
\def\,{\thinspace}\thereflist\par}}
%
%   Next a footnote macro (with automatic numbering) which sets the
%   footnote  \small.
%
%   Typical use :
%         ..... are yellow.\fnote{By yellow here we mean Britsh
%    Standard colour BS3320.} 
%
\newcount\footnotenumber         % Register for footnote number
\footnotenumber=1                % set initially to 1
\def\fnote#1{\xdef\nextkey{\number\footnotenumber}%
{\small\ifnum\footnotenumber>9\parindent=14pt%
\else\parindent=10pt\fi\footnote{$^{\number\footnotenumber}$}%
{\hglue-5pt#1}\global\advance\footnotenumber by 1}}
%
%
%   Next macros for handling figures with automatic numbering (using 
%   TeX's \midinsert to float the figure to a suitable place).
%   
%   The \figure ... \endfigure macro centres the figure and adds
%   an automatically numbered label  Figure XX  after it.
%
%   If you have a caption, then type \caption{caption text} 
%   somewhere between \figure and \endfigure.  The macro
%   will then add  Figure XX: caption text  after the figure.
%
%   If you want an unnumbered or uncentred figure, then use TeX's raw 
%       \midinsert Figure instructions \endinsert  
%   and if you want a numbered figure label in the same style then
%   use \caption{caption text} outside of  \figure ... \endfigure.
%
%   If you need just the label Figure XX  outside of  \figure ... \endfigure
%   then type  \figurelabel .
%
\newcount\figurenumber          % register for figure number
\figurenumber=1                 % set initially to 1
\def\caption#1{\xdef\nextkey{\number\figurenumber}%
\cl{\small Figure \number\figurenumber: #1}%
\global\advance\figurenumber by 1}
\def\figurelabel{\xdef\nextkey{\number\figurenumber}%
\cl{\small Figure \number\figurenumber}%
\global\advance\figurenumber by 1}
\long\def\figure#1\endfigure{{\xdef\nextkey{\number\figurenumber}%
\let\captiontext\relax\def\caption##1{\xdef\captiontext{##1}}%
\midinsert\cl{\ignorespaces#1\unskip\unskip\unskip\unskip}\vglue6pt\cl{\small 
Figure \number\figurenumber\ifx\captiontext\relax\else: \captiontext
\fi}\endinsert\global\advance\figurenumber by 1}}
%
%   Macros for self-correcting internal references.
%
%   There are two macros  \key{KEY}  and  \ref{KEY} .
%
%   The \key macro sets up KEY as a key for whatever number is 
%   being referenced and the \ref macro converts the KEY into 
%   that number.  Type \key after a  \section or \proc or 
%   \label or \fnote or \figure or \caption or \figurelabel .
%
%   Example:
%
%       \section{Introduction}\key{intro}
%       \proc{Theorem}\key{MainTh}Lemons are yelloy\endproc
%       Here we follow\fnote{Follow in the sense of Dickens}
%       \key{Dickens-note}the crowd ....  
%
%       In section \ref{intro}
%       we stated theorem \ref{mainTh} and noted (see footnote 
%       \ref{Dickens-note}) ...
%
\def\nextkey{??}   %  initialise \nextkey (which is reset by all the
%                     numbering macros)
%
\def\key#1{\expandafter\xdef\csname #1\endcsname{\nextkey}}
\def\ref#1{\expandafter\ifx\csname #1\endcsname\relax
\immediate\write16{Reference {#1} undefined}??\else
\csname #1\endcsname\fi}
%
%   Note:  If the KEY contains only letters then \KEY has exactly the
%   same meaning as \ref{KEY} so in the example you could have:
%
%       In section \intro\ we ....
%
%   The \key will work at any time after the macro which sets the
%   number, provided no other macro which sets a number has been used. 
%
%   Macros for forward references:
%              =======
%   The \key \ref macros ONLY work for backwards references.  If you  
%   want to use forwards references, then type \useforwardrefs  near
%   the beginning of your file.  The KEY's are then stored in an
%   auxiliary  .ref  file and you then suffer the same disadvantage as
%   when using LaTeX that you must TeX the file twice to get
%   the references correct.
%
%   To use a forward ref type \ref{KEY}.  (You can type the
%   alternative  \KEY  but you'll get an error on first TeX'ing 
%   if the \KEY is not yet defined.) 
%
%   The macro also allows external references to be listed at the end 
%   of the file (if you wish to).  (Indeed they can be typed anywhere
%   before the \references command.)  You can combine the reference list
%   and the \references command by typing the references (using the
%   same syntax as before) between the commands \biblio and \endbiblio 
%   (don't type \references or they'll be printed twice).
%
\newread\gtinfile
\newwrite\gtreffile
\def\useforwardrefs{
\openin\gtinfile\jobname.ref
\ifeof\gtinfile
\closein\gtinfile
\immediate\write16{No file \jobname.ref}
\else
\closein\gtinfile
\input \jobname.ref
\fi
\immediate\openout\gtreffile \jobname.ref
%
%   Adapt \key :
%
\def\key##1{{\def\\{\noexpand}%
\expandafter\xdef\csname ##1\endcsname{\nextkey}%
\immediate\write\gtreffile{\\\expandafter\\\def\\\csname ##1\\\endcsname%
{\nextkey}}}}
%
%  Adapt macros for external references:  
%
\long\def\reflist##1\endreflist{%
\long\def\thereflist{##1}{\def\refkey####1####2\par{\xdef####1{%
\number\refnumber}{\def\\{\noexpand}\immediate\write\gtreffile
{\\\def\\####1{\number\refnumber}}}\global\advance\refnumber by 1}%
\def\key####1####2\par{\expandafter\xdef%
\csname####1\endcsname{\number\refnumber}%
{\def\\{\noexpand}\immediate\write\gtreffile
{\\\expandafter\\\def\\\csname ####1\\\endcsname{\number\refnumber}}}
\global\advance\refnumber by 1}##1\par}}
\long\def\biblio##1\endbiblio{\reflist##1\endreflist\references}%
%
%  Adapt obselete key macros (\numkey, \seckey and \figkey):
%
\def\numkey##1{{\def\\{\noexpand}%
\xdef##1{\number\sectionnumber.\number\resultnumber}
\immediate\write\gtreffile{\\\def\\##1%
{\number\sectionnumber.\number\resultnumber}}}}
\def\seckey##1{{\def\\{\noexpand}\xdef##1{\number\sectionnumber}
\immediate\write\gtreffile{\\\def\\##1{\number\sectionnumber}}}}
\def\figkey##1{\xdef##1{\number\figurenumber}%
{\def\\{\noexpand}\immediate\write\gtreffile%
{\\\def\\##1{\number\figurenumber}}}
\number\figurenumber\global\advance\figurenumber by 1}
}   %  end of \useforwardrefs
%
%
%   The next five macros are obselete and have been superseeded by
%   the general \key macro above.  They are included merely to 
%   maintain backward compatibility for the package:
%
%
\def\figkey#1{\xdef#1{\number\figurenumber}%
\number\figurenumber\global\advance\figurenumber by 1}
\def\fig#1#2\endfig{%
\midinsert\cl{#2}\vglue6pt\cl{\small Figure #1}\endinsert}
\def\newfig{\number\figurenumber\global\advance\figurenumber by 1}
\def\numkey#1{\xdef#1{\number\sectionnumber.\number\resultnumber}}
\def\seckey#1{\xdef#1{\number\sectionnumber}}
%
%   End of obselete macros.
%
%
%   The next macro is a version of the verbatim macro given by Knuth.
%
%   This macro produces a "verbatim" printout of
%   any ASCII string which does not contain the symbol "
%   (TeX files do not usually contain " 's).
%   More precisely, everything between consecutive pairs
%   of " 's is printed verbatim in the typewriter font cmtt.
%   For an explanation of how the macro works, see Knuth pp 420-1.
%
%   There are two switches: \verb (which switches the macro on)
%   and \brev which switches the macro off (the default).  When
%   the macro is switched off the symbol " has its usual 
%   meaning for TeX.  To use the macro, type \verb before use
%   and the use " to switch verbatim on and off.  Be careful
%   not to use " for any other purpose.  There is no need to
%   switch the macro off again unless you need to use " for
%   some other purpose (eg making  \mathchardef 's).  Note 
%   that the macro MUST BE OFF before inputting  amsnames.tex .
%
%   Whether the macro is on or off you can always use the
%   control sequence \dq (double quote) for " e.g.
%   \mathchardef\sum=\dq1350  is perfectly valid.
%   The control sequence \ttq is an abbreviation for
%   {\tt\dq}.  Thus "\ttq" will produce " (in cmtt)
%   inside a verbatim quote.
%
%
   %  define a code for " so it can be used when \verb is on
  %  code for " in cmtt
%
\def\verb{\catcode`\"=\active}       %  The main
\def\brev{\catcode`\"=12}            %  switches.
\brev                                %  Prime switches and
\verb                                %  switch on.
{\obeyspaces\gdef {\ }}              
{\catcode`\`=\active\gdef`{\relax\lq}}
\def"{%
\begingroup\baselineskip=12pt\def\par{\leavevmode\endgraf}%
\tt\obeylines\obeyspaces\parskip=0pt\parindent=0pt%
\catcode`\$=12\catcode`\&=12\catcode`\^=12\catcode`\#=12%
\catcode`\_=12\catcode`\~=12%
\catcode`\{=12\catcode`\}=12\catcode`\%=12\catcode`\\=12%
\catcode`\`=\active\let"\endgroup}
\brev      %   Finally switch the macro off (for safety)
%
%   Macros for itemised lists.   Typical use :
%    
%    \items
%    \item{(i)}Colours must be defined.
%    \item{(ii)}Colour cards may not be cited.
%    \enditems
%
%   Result :
%
%    (i)  Colours must be defined. 
%   (ii)  Colour cards may not be cited.
%
%
\def\items{\par\leftskip = 25pt}           % Start of itemised list         
\def\enditems{\par\leftskip = 0pt}         % end of itemised list   
\def\item#1{\par\leavevmode\llap{#1\stdspace}%
\ignorespaces}                             % labelled item
\def\itemb{\item{$\bullet$}}               % bulleted item.
%
%   The \quote ... \endquote macros are for typesetting quotations :
%

%
%   A few useful abbreviations :
%
\def\co{\colon\thinspace}    %  Colon with correct spacing for maps.
\def\np{\vfil\eject}         %  Forced page break (new page).
\def\nl{\hfil\break}         %  New line.
\def\cl{\centerline}         %  Centerline
        %  The journal title in recommended style
    %  for monographs
\def\agt{{\mathsurround=0pt\it$\cal A\mskip-.7mu$lgebraic \&\ 
$\cal G\mskip-2mu$eometric $\cal T\!\!$opology}}  % AGT
%
%    Finally some macros for automatic title page or header generation.
%    To use them type your header information using the following  
%    example as a guide :
%
%    Note that \\ is used as standard separator (for lines in \title and
%    \address, between authors and between email addresses or URL's)
%    and that \email, \url and \secondaddress are optional.
%

% Example:  \title{A short spoof paper\\with a two-line title}
% =======   \authors{Albert Einstein\\Leonardo da Vinci}
%           \address{IAS\\Princeton}\secondaddress{Renaissance\\Venice}
%           \email{ae@ias.princeton.edu\\ldv@ren.ven.hist}
%           \abstract 
%           A short spoof paper with a very short abstract.
%           \endabstract 
%           \primaryclass{00-01, 00-02}\secondaryclass{68-00, 68-01}
%           \keywords{Short, spoof, paper}
%           \maketitlepage
%
%
%    The title page or header will then be generated automatically.
%
%
%    Define the various ingredients of the title page:
%
\def\title#1{\def\thetitle{#1}}

\def\author#1{\edef\previousauthors{\theauthors}
 \ifx\theauthors\relax\def\theauthors{#1}\else
 \def\theauthors{\previousauthors\par#1}\fi}

\let\authors\author        % aliases
\def\address#1{\edef\previousaddresses{\theaddress}
 \ifx\theaddress\relax\def\theaddress{#1}\else
 \def\theaddress{\previousaddresses\par\vskip 2pt\par#1}\fi}
                             % alias
\def\secondaddress#1{\edef\previousaddresses{\theaddress}
 \ifx\theaddress\relax\def\theaddress{#1}\else
 \def\theaddress{\previousaddresses\par{\rm and}\par#1}\fi}   

\def\email#1{\edef\previousemails{\theemail}
 \ifx\theemail\relax\def\theemail{#1}\else
 \def\theemail{\previousemails\hskip 0.75em\relax#1}\fi}
  % aliases
\def\secondemail#1{\edef\previousemails{\theemail}
 \ifx\theemail\relax\def\theemail{#1}\else
 \def\theemail{\previousemails\hskip 0.75em{\rm and}\hskip 0.75em
 \relax#1}\fi}
\def\url#1{\edef\previousurls{\theurl}
 \ifx\theurl\relax\def\theurl{#1}\else
 \def\theurl{\previousurls\hskip 0.75em\relax#1}\fi}
      % aliases
\def\secondurl#1{\edef\previousurls{\theurl}
 \ifx\theurl\relax\def\theurl{#1}\else
 \def\theurl{\previousurls\hskip 0.75em{\rm and}\hskip 0.75em
 \relax#1}\fi}
\long\def\abstract#1\endabstract{\long\def\theabstract{#1}}
\def\primaryclass#1{\def\theprimaryclass{#1}}
                        % alias
\def\secondaryclass#1{\def\thesecondaryclass{#1}}
\def\keywords#1{\def\thekeywords{#1}}
%
%  Set \\ to \par and title page items to \relax to initialise macros :
%
\let\\\par\let\thetitle\relax\let\theshorttitle\relax
\let\theauthors\relax\let\theshortauthors\relax
\let\theaddress\relax\let\theshortaddress\relax
\let\theemail\relax\let\theurl\relax
\let\theabstract\relax\let\theprimaryclass\relax
\let\thesecondaryclass\relax\let\thekeywords\relax
%
%
%
%   Basic title page layout (edit this macro if you
%   wish to adjust the title page layout) :
%
\long\def\maketitlepage{    % start of definition of \maketitlepage

\vglue 0.2truein   % top margin

% title :
%
{\parskip=0pt\leftskip 0pt plus 1fil\def\\{\par\smallskip}{\large
\bf\thetitle}\par\medskip}   

\vglue 0.15truein 

% authors :
%
{\parskip=0pt\leftskip 0pt plus 1fil\def\\{\par}{\sc\theauthors}
\par\medskip}%
 
\vglue 0.1truein 

% address(es) email's and URL's (with switches to detect whether the
% optional items have been used) :
%
{\small\parskip=0pt
{\leftskip 0pt plus 1fil\def\\{\par}{\sl\theaddress}\par}
\ifx\theemail\relax\else  % email address?
\vglue 5pt \def\\{\stdspace{\rm and}\stdspace} 
\cl{Email:\stdspace\tt\theemail}\fi
\ifx\theurl\relax\else    % URL given?
\vglue 5pt \def\\{\stdspace{\rm and}\stdspace} 
\cl{URL:\stdspace\tt\theurl}\fi\par}

\vglue 7pt 

{\bf Abstract}

\vglue 5pt

\theabstract

\vglue 7pt 

{\bf AMS Classification numbers}\quad Primary:\quad \theprimaryclass\par

Secondary:\quad \thesecondaryclass

\vglue 5pt 

{\bf Keywords:}\quad \thekeywords

\np  % page break at the end of the title page

}    % end of definition of \maketitlepage
%
%    % \makeshorttitle (for general preprints) doesn't take a new page
%
\long\def\makeshorttitle{    % start of definition of \makeshorttitle

%\vglue 0.2truein   % top margin

% title :
%
{\parskip=0pt\leftskip 0pt plus 1fil\def\\{\par\smallskip}{\large
\bf\thetitle}\par\medskip}   

\vglue 0.05truein 

% authors :
%
{\parskip=0pt\leftskip 0pt plus 1fil\def\\{\par}{\sc\theauthors}
\par\medskip}%
 
\vglue 0.03truein 

% address(es) email's and URL's (with switches to detect whether the
% optional items have been used) :
%
{\small\parskip=0pt
{\leftskip 0pt plus 1fil\def\\{\par}{\sl\ifx\theshortaddress\relax
\theaddress\else\theshortaddress\fi}\par}
\ifx\theemail\relax\else  % email address?
\vglue 5pt \def\\{\stdspace{\rm and}\stdspace} 
\cl{Email:\stdspace\tt\theemail}\fi
\ifx\theurl\relax\else    % URL given?
\vglue 5pt \def\\{\stdspace{\rm and}\stdspace} 
\cl{URL:\stdspace\tt\theurl}\fi\par}

\vglue 10pt 

% abstract and classification numbers (with switches):

{\small\leftskip 25pt\rightskip 25pt{\bf Abstract}\stdspace\theabstract

{\bf AMS Classification}\stdspace\theprimaryclass
\ifx\thesecondaryclass\relax\else; \thesecondaryclass\fi\par
{\bf Keywords}\stdspace \thekeywords\par}
\vglue 7pt
}    % end of definition of \makeshorttitle
\let\maketitle\makeshorttitle        %% alias
%
%    %%%% \makeagttitle (for AGT) similar to \makeshorttitle but
%         with addresses omitted (they go at the end)
%
%%%% publication info and test defaults:

\def\volumenumber#1{\def\thevolumenumber{#1}}
\def\volumeyear#1{\def\thevolumeyear{#1}}
\def\pagenumbers#1#2{\def\startpage{#1}\def\finishpage{#2}}
\def\published#1{\def\publishdate{#1}}
\def\received#1{\def\receiveddate{#1}}
\def\revised#1{\def\reviseddate{#1}}
\let\reviseddate\relax
%% Defaults for authors to use to check layout
\volumenumber{X}
\volumeyear{20XX}
\pagenumbers{1}{XXX}
\published{XX Xxxember 20XX}

\long\def\makeagttitle{   %%% start of definition of \makeagttitle
\agt\hfill      %   Journal title (top left) 
%   logo placeholder (top right)
\hbox to 60truept{\vbox to 0pt{\vglue -14truept{\bf [Logo here]}\vss}\hss}
\break
{\small Volume \thevolumenumber\ (\thevolumeyear)
\startpage--\finishpage\nl
Published: \publishdate}

\vglue .2truein

% title
{\parskip=0pt\leftskip 0pt plus 1fil\def\\{\par\smallskip}{\large
\bf\thetitle}\par\medskip}   
\vglue 0.05truein 

% authors :
%
{\parskip=0pt\leftskip 0pt plus 1fil\def\\{\par}{\sc\theauthors}
\par\medskip}%
 
\vglue 0.03truein 

%  abstract and classification numbers:

{\small\leftskip 25truept\rightskip 25truept{\bf Abstract}\stdspace\theabstract

{\bf AMS Classification}\stdspace\theprimaryclass
\ifx\thesecondaryclass\relax\else; \thesecondaryclass\fi\par
{\bf Keywords}\stdspace \thekeywords\par}\vglue 7truept

}   %%%% end of definition of \makeagttitle

%%%%% Macro to typeset addresses (typically at the end of the paper)

\def\Addresses{\bigskip
{\small \parskip 0pt \leftskip 0pt \rightskip 0pt plus 1fil \def\\{\par}
\sl\theaddress\par\medskip \rm Email:\stdspace\tt\theemail\par
\ifx\theurl\relax\else\smallskip \rm URL:\stdspace\tt\theurl\par\fi}}

\def\agtart{%   Full mock-up of AGT article style (for authors to test with)
%  get print centerpage:
\hoffset 14truemm
\voffset 31truemm
\font\phead=cmsl9 scaled 950
\font\pnum=cmbx10 scaled 913
\font\pfoot=cmsl9 scaled 950
%  headline and footline
\headline{\vbox to 0pt{\vskip -4.5mm\line{\small\phead\ifnum
\count0=\startpage ISSN numbers are printed here
\hfill {\pnum\folio}\else\ifodd\count0\def\\{ }% 
\ifx\theshorttitle\relax\thetitle\else\theshorttitle\fi\hfill{\pnum\folio}
\else\def\\{ and }{\pnum\folio}\hfill\ifx\theshortauthors\relax\theauthors
\else\theshortauthors\fi\fi\fi}\vss}}
\footline{\vbox to 0pt{\vglue 0mm\line{\small\pfoot\ifnum\count0=\startpage
Copyright declaration is printed here\hfill\else
\agt, Volume \thevolumenumber\ (\thevolumeyear)\hfill\fi}\vss}}
%  force \agttitle
\let\maketitle\makeagttitle\let\makeshorttitle\makeagttitle}

%%%
%%%  This is agtout.tex.  
%%%
%%%  This the version of  gtoutput.tex  intended to finish formatting
%%%  papers published in Algebriac & Geometric Topology and stored in the
%%%  arXiv.   All versions of  gtoutput.tex  are copyright 
%%%  GT Publications and are to be used _only_ for formatting
%%%  the officially published version of ABT or G&T papers.
%%%
%%%
%%%                                             Colin Rourke  27.102000
%%%
%%%  To create header file  head.xxx  comment out the first \endinput

%  test for latex or plain tex
\def\ifplaintex{\expandafter\ifx\csname documentclass\endcsname\relax}

\def\gtp{{\mathsurround=0pt\it $\cal G\mskip-2mu$eometry \&\ 
$\cal T\!\!$opology $\cal P\!$ublications}}  % GT publications

\def\recd{{\small Received:\qua\receiveddate\ifx\reviseddate\relax
\else\qquad Revised:\qua\reviseddate\fi\par}} 

%  define the various new ingredients of the title page and the data
%  output files

\def\lognumber#1{\def\thelognumber{#1}}
\def\volumenumber#1{\def\thevolumenumber{#1}}
\def\volumeyear#1{\def\thevolumeyear{#1}}
\def\papernumber#1{\def\thepapernumber{#1}}
\def\pagenumbers#1#2{\def\startpage{#1}\def\finishpage{#2}}
\def\published#1{\def\publishdate{#1}}

\def\received#1{\def\receiveddate{#1}}
\def\revised#1{\def\reviseddate{#1}}
\def\accepted#1{\def\accepteddate{#1}}
\def\asciititle#1{\def\theasciititle{#1}}

\def\asciiaddress#1{\def\theasciiaddress{#1}}

\long\def\asciiabstract#1{\long\def\theasciiabstract{#1}}

%  initialise

\let\\\par\let\thelognumber\relax\let\thevolumenumber\relax
\let\thepapernumber\relax\let\thevolumeyear\relax\let\startpage\relax
\let\finishpage\relax\let\publishdate\relax\let\receiveddate\relax
\let\reviseddate\relax\let\accepteddate\relax\let\theasciititle\relax
\let\theasciiauthors\relax\let\theasciiaddress\relax
\let\theasciiabstract\relax

\let\theasciiemail\relax

%%%% fonts for AGT logo:

\ifplaintex
\font\logobig=cmssbx10 scaled 3836
\font\logomed=cmssbx10 scaled 2557
\else
\font\logobig=cmssbx10 scaled 4200
\font\logomed=cmssbx10 scaled 2800
\fi

\long\def\makeagttitle{   %%% start of definition of \makeagttitle
\count0=\startpage
\agt\hfill      %   Journal title (top left) 
%   logo (top right)
\hbox to 45truept{\vbox to 0pt{\vglue -13truept{\logomed A\kern -.37em{\logobig 
T}\kern -.38em G}\vss}\hss}
\break
{\small Volume \thevolumenumber\ (\thevolumeyear)
\startpage--\finishpage\nl
Published: \publishdate}

\vglue .25truein

% title
{\parskip=0pt\leftskip 0pt plus
1fil\def\\{\par\smallskip}{\Large\bf\thetitle}\par\medskip} \vglue
0.05truein

% authors :
%
{\parskip=0pt\leftskip 0pt plus 1fil\def\\{\par}{\sc\theauthors}
\par\medskip}%
 
\vglue 0.03truein 

%  abstract and classification numbers:

{\small\leftskip 25truept\rightskip 25truept{\bf Abstract}\stdspace\theabstract

{\bf AMS Classification}\stdspace\theprimaryclass
\ifx\thesecondaryclass\relax\else; \thesecondaryclass\fi\par
{\bf Keywords}\stdspace \thekeywords\par}\vglue 7truept

}   %%%% end of definition of \makeagttitle

\ifplaintex
%  get print centerpage:
\hoffset 14truemm
\voffset 31truemm
%  fonts for headline and footline
\font\phead=cmsl9 scaled 950
\font\pnum=cmbx10 scaled 913
\font\pfoot=cmsl9 scaled 950
%  headline and footline
\headline{\vbox to 0pt{\vskip -4.5mm\line{\small\phead\ifnum
\count0=\startpage ISSN 1472-2739 (on-line) 1472-2747 (printed)
\hfill {\pnum\folio}\else\ifodd\count0\def\\{ }% 
\ifx\theshorttitle\relax\thetitle\else\theshorttitle\fi\hfill{\pnum\folio}
\else\def\\{ and }{\pnum\folio}\hfill\ifx\theshortauthors\relax\theauthors
\else\theshortauthors\fi\fi\fi}\vss}}
\footline{\vbox to 0pt{\vglue 0mm\line{\small\pfoot\ifnum\count0=\startpage
\copyright\ \gtp\hfill\else
\agt, Volume \thevolumenumber\ (\thevolumeyear)\hfill\fi}\vss}}
\else
%  get print centerpage:
\headsep 23pt
\footskip 35pt
\hoffset -4truemm
\voffset 12.5truemm
%  fonts for headline and footline
\font\lhead=cmsl9 scaled 1050
\font\lnum=cmbx10 
\font\lfoot=cmsl9 scaled 1050
\makeatletter
%  headline and footline
\def\@oddhead{{\small\lhead\ifnum\count0=\startpage ISSN 1472-2739 
(on-line) 1472-2747 (printed)\hfill {\lnum\number\count0}\else\ifodd\count0
\def\\{ }\ifx\theshorttitle\relax \thetitle \else\theshorttitle\fi\hfill
{\lnum\number\count0}\else\def\\{ and }{\lnum\number\count0}
\hfill\ifx\theshortauthors\relax 
\theauthors\else\theshortauthors\fi\fi\fi}}\def\@evenhead{\@oddhead}
\def\@oddfoot{\small\lfoot\ifnum\count0=\startpage\copyright\ \gtp\hfill\else
\agt, Volume \thevolumenumber\ (\thevolumeyear)\hfill\fi}
\def\@evenfoot{\@oddfoot}
\makeatother
\fi
%  force \makeagttitle
\let\maketitlepage\makeagttitle
\let\makeshorttitle\maketitlepage
\let\maketitle\maketitlepage

   %%%comment out to create xxx header file

\newwrite\gtoutfile
\long\gdef\makeheadfile{  %%% start of definition of \makeheadfile
{\def\\{, }\def\s{ }
\immediate\openout\gtoutfile head.xxx
\immediate\write\gtoutfile{To: math@arxiv.org}
\immediate\write\gtoutfile{Subject: put OR rep NNNNN:ppppp}
\immediate\write\gtoutfile{--text follows this line--}
\immediate\write\gtoutfile{Proxy-for: \ifx\theasciiauthors\relax
\theauthors\else\theasciiauthors\fi\s<\ifx\theasciiemail\relax\theemail\else\theasciiemail\fi>}
\immediate\write\gtoutfile{\noexpand\\}
\immediate\write\gtoutfile{Authors: \ifx\theasciiauthors\relax
\theauthors\else\theasciiauthors\fi}
{\def\\{ }\immediate\write\gtoutfile{Title: \ifx\theasciititle\relax
\thetitle\else\theasciititle\fi}}
\immediate\write\gtoutfile{Subj-class: GT or SG, GR etc}
\immediate\write\gtoutfile{MSC-class: \theprimaryclass\ifx\thesecondaryclass\relax\else, \thesecondaryclass\fi}
\immediate\write\gtoutfile{Journal-ref: Algebr. Geom. Topol. \thevolumenumber\s
(\thevolumeyear) \startpage-\finishpage}
\immediate\write\gtoutfile{Comments: Published by Algebraic and
Geometric Topology at}
\immediate\write\gtoutfile{\s\s\s  http://www.maths.warwick.ac.uk/agt/AGTVol\thevolumenumber/agt-\thevolumenumber-\thepapernumber.abs.html}
\immediate\write\gtoutfile{\noexpand\\}
\immediate\write\gtoutfile{}
\ifx\theasciiabstract\relax
\immediate\write\gtoutfile{\theabstract}\else
\immediate\write\gtoutfile{\theasciiabstract}\fi
\immediate\write\gtoutfile{}
\immediate\write\gtoutfile{\noexpand\\}
\immediate\write\gtoutfile{}
\immediate\closeout\gtoutfile}}  %%% end of definition of \makeheadfile

\def\maketitlepage{\makeagttitle\makeheadfile}
\let\makeshorttitle\maketitlepage
\let\maketitle\maketitlepage

%%%
%%%  This is agtout.tex.  
%%%
%%%  This the version of  gtoutput.tex  intended to finish formatting
%%%  papers published in Algebriac & Geometric Topology and stored in the
%%%  arXiv.   All versions of  gtoutput.tex  are copyright 
%%%  GT Publications and are to be used _only_ for formatting
%%%  the officially published version of ABT or G&T papers.
%%%
%%%
%%%                                             Colin Rourke  27.102000
%%%
%%%  To create header file  head.xxx  comment out the first \endinput

%  test for latex or plain tex
\def\ifplaintex{\expandafter\ifx\csname documentclass\endcsname\relax}

\def\gtp{{\mathsurround=0pt\it $\cal G\mskip-2mu$eometry \&\ 
$\cal T\!\!$opology $\cal P\!$ublications}}  % GT publications

\def\recd{{\small Received:\qua\receiveddate\ifx\reviseddate\relax
\else\qquad Revised:\qua\reviseddate\fi\par}} 

%  define the various new ingredients of the title page and the data
%  output files

\def\lognumber#1{\def\thelognumber{#1}}
\def\volumenumber#1{\def\thevolumenumber{#1}}
\def\volumeyear#1{\def\thevolumeyear{#1}}
\def\papernumber#1{\def\thepapernumber{#1}}
\def\pagenumbers#1#2{\def\startpage{#1}\def\finishpage{#2}}
\def\published#1{\def\publishdate{#1}}

\def\received#1{\def\receiveddate{#1}}
\def\revised#1{\def\reviseddate{#1}}
\def\accepted#1{\def\accepteddate{#1}}
\def\asciititle#1{\def\theasciititle{#1}}

\def\asciiaddress#1{\def\theasciiaddress{#1}}

\long\def\asciiabstract#1{\long\def\theasciiabstract{#1}}

%  initialise

\let\\\par\let\thelognumber\relax\let\thevolumenumber\relax
\let\thepapernumber\relax\let\thevolumeyear\relax\let\startpage\relax
\let\finishpage\relax\let\publishdate\relax\let\receiveddate\relax
\let\reviseddate\relax\let\accepteddate\relax\let\theasciititle\relax
\let\theasciiauthors\relax\let\theasciiaddress\relax
\let\theasciiabstract\relax

\let\theasciiemail\relax

%%%% fonts for AGT logo:

\ifplaintex
\font\logobig=cmssbx10 scaled 3836
\font\logomed=cmssbx10 scaled 2557
\else
\font\logobig=cmssbx10 scaled 4200
\font\logomed=cmssbx10 scaled 2800
\fi

\long\def\makeagttitle{   %%% start of definition of \makeagttitle
\count0=\startpage
\agt\hfill      %   Journal title (top left) 
%   logo (top right)
\hbox to 45truept{\vbox to 0pt{\vglue -13truept{\logomed A\kern -.37em{\logobig 
T}\kern -.38em G}\vss}\hss}
\break
{\small Volume \thevolumenumber\ (\thevolumeyear)
\startpage--\finishpage\nl
Published: \publishdate}

\vglue .25truein

% title
{\parskip=0pt\leftskip 0pt plus
1fil\def\\{\par\smallskip}{\Large\bf\thetitle}\par\medskip} \vglue
0.05truein

% authors :
%
{\parskip=0pt\leftskip 0pt plus 1fil\def\\{\par}{\sc\theauthors}
\par\medskip}%
 
\vglue 0.03truein 

%  abstract and classification numbers:

{\small\leftskip 25truept\rightskip 25truept{\bf Abstract}\stdspace\theabstract

{\bf AMS Classification}\stdspace\theprimaryclass
\ifx\thesecondaryclass\relax\else; \thesecondaryclass\fi\par
{\bf Keywords}\stdspace \thekeywords\par}\vglue 7truept

}   %%%% end of definition of \makeagttitle

\ifplaintex
%  get print centerpage:
\hoffset 14truemm
\voffset 31truemm
%  fonts for headline and footline
\font\phead=cmsl9 scaled 950
\font\pnum=cmbx10 scaled 913
\font\pfoot=cmsl9 scaled 950
%  headline and footline
\headline{\vbox to 0pt{\vskip -4.5mm\line{\small\phead\ifnum
\count0=\startpage ISSN 1472-2739 (on-line) 1472-2747 (printed)
\hfill {\pnum\folio}\else\ifodd\count0\def\\{ }% 
\ifx\theshorttitle\relax\thetitle\else\theshorttitle\fi\hfill{\pnum\folio}
\else\def\\{ and }{\pnum\folio}\hfill\ifx\theshortauthors\relax\theauthors
\else\theshortauthors\fi\fi\fi}\vss}}
\footline{\vbox to 0pt{\vglue 0mm\line{\small\pfoot\ifnum\count0=\startpage
\copyright\ \gtp\hfill\else
\agt, Volume \thevolumenumber\ (\thevolumeyear)\hfill\fi}\vss}}
\else
%  get print centerpage:
\headsep 23pt
\footskip 35pt
\hoffset -4truemm
\voffset 12.5truemm
%  fonts for headline and footline
\font\lhead=cmsl9 scaled 1050
\font\lnum=cmbx10 
\font\lfoot=cmsl9 scaled 1050
\makeatletter
%  headline and footline
\def\@oddhead{{\small\lhead\ifnum\count0=\startpage ISSN 1472-2739 
(on-line) 1472-2747 (printed)\hfill {\lnum\number\count0}\else\ifodd\count0
\def\\{ }\ifx\theshorttitle\relax \thetitle \else\theshorttitle\fi\hfill
{\lnum\number\count0}\else\def\\{ and }{\lnum\number\count0}
\hfill\ifx\theshortauthors\relax 
\theauthors\else\theshortauthors\fi\fi\fi}}\def\@evenhead{\@oddhead}
\def\@oddfoot{\small\lfoot\ifnum\count0=\startpage\copyright\ \gtp\hfill\else
\agt, Volume \thevolumenumber\ (\thevolumeyear)\hfill\fi}
\def\@evenfoot{\@oddfoot}
\makeatother
\fi
%  force \makeagttitle
\let\maketitlepage\makeagttitle
\let\makeshorttitle\maketitlepage
\let\maketitle\maketitlepage

   %%%comment out to create xxx header file

\newwrite\gtoutfile
\long\gdef\makeheadfile{  %%% start of definition of \makeheadfile
{\def\\{, }\def\s{ }
\immediate\openout\gtoutfile head.xxx
\immediate\write\gtoutfile{To: math@arxiv.org}
\immediate\write\gtoutfile{Subject: put OR rep NNNNN:ppppp}
\immediate\write\gtoutfile{--text follows this line--}
\immediate\write\gtoutfile{Proxy-for: \ifx\theasciiauthors\relax
\theauthors\else\theasciiauthors\fi\s<\ifx\theasciiemail\relax\theemail\else\theasciiemail\fi>}
\immediate\write\gtoutfile{\noexpand\\}
\immediate\write\gtoutfile{Authors: \ifx\theasciiauthors\relax
\theauthors\else\theasciiauthors\fi}
{\def\\{ }\immediate\write\gtoutfile{Title: \ifx\theasciititle\relax
\thetitle\else\theasciititle\fi}}
\immediate\write\gtoutfile{Subj-class: GT or SG, GR etc}
\immediate\write\gtoutfile{MSC-class: \theprimaryclass\ifx\thesecondaryclass\relax\else, \thesecondaryclass\fi}
\immediate\write\gtoutfile{Journal-ref: Algebr. Geom. Topol. \thevolumenumber\s
(\thevolumeyear) \startpage-\finishpage}
\immediate\write\gtoutfile{Comments: Published by Algebraic and
Geometric Topology at}
\immediate\write\gtoutfile{\s\s\s  http://www.maths.warwick.ac.uk/agt/AGTVol\thevolumenumber/agt-\thevolumenumber-\thepapernumber.abs.html}
\immediate\write\gtoutfile{\noexpand\\}
\immediate\write\gtoutfile{}
\ifx\theasciiabstract\relax
\immediate\write\gtoutfile{\theabstract}\else
\immediate\write\gtoutfile{\theasciiabstract}\fi
\immediate\write\gtoutfile{}
\immediate\write\gtoutfile{\noexpand\\}
\immediate\write\gtoutfile{}
\immediate\closeout\gtoutfile}}  %%% end of definition of \makeheadfile

\def\maketitlepage{\makeagttitle\makeheadfile}
\let\makeshorttitle\maketitlepage
\let\maketitle\maketitlepage

%%%
%%%  This is agtout.tex.  
%%%
%%%  This the version of  gtoutput.tex  intended to finish formatting
%%%  papers published in Algebriac & Geometric Topology and stored in the
%%%  arXiv.   All versions of  gtoutput.tex  are copyright 
%%%  GT Publications and are to be used _only_ for formatting
%%%  the officially published version of ABT or G&T papers.
%%%
%%%
%%%                                             Colin Rourke  27.102000
%%%
%%%  To create header file  head.xxx  comment out the first \endinput

%  test for latex or plain tex
\def\ifplaintex{\expandafter\ifx\csname documentclass\endcsname\relax}

\def\gtp{{\mathsurround=0pt\it $\cal G\mskip-2mu$eometry \&\ 
$\cal T\!\!$opology $\cal P\!$ublications}}  % GT publications

\def\recd{{\small Received:\qua\receiveddate\ifx\reviseddate\relax
\else\qquad Revised:\qua\reviseddate\fi\par}} 

%  define the various new ingredients of the title page and the data
%  output files

\def\lognumber#1{\def\thelognumber{#1}}
\def\volumenumber#1{\def\thevolumenumber{#1}}
\def\volumeyear#1{\def\thevolumeyear{#1}}
\def\papernumber#1{\def\thepapernumber{#1}}
\def\pagenumbers#1#2{\def\startpage{#1}\def\finishpage{#2}}
\def\published#1{\def\publishdate{#1}}

\def\received#1{\def\receiveddate{#1}}
\def\revised#1{\def\reviseddate{#1}}
\def\accepted#1{\def\accepteddate{#1}}
\def\asciititle#1{\def\theasciititle{#1}}

\def\asciiaddress#1{\def\theasciiaddress{#1}}

\long\def\asciiabstract#1{\long\def\theasciiabstract{#1}}

%  initialise

\let\\\par\let\thelognumber\relax\let\thevolumenumber\relax
\let\thepapernumber\relax\let\thevolumeyear\relax\let\startpage\relax
\let\finishpage\relax\let\publishdate\relax\let\receiveddate\relax
\let\reviseddate\relax\let\accepteddate\relax\let\theasciititle\relax
\let\theasciiauthors\relax\let\theasciiaddress\relax
\let\theasciiabstract\relax

\let\theasciiemail\relax

%%%% fonts for AGT logo:

\ifplaintex
\font\logobig=cmssbx10 scaled 3836
\font\logomed=cmssbx10 scaled 2557
\else
\font\logobig=cmssbx10 scaled 4200
\font\logomed=cmssbx10 scaled 2800
\fi

\long\def\makeagttitle{   %%% start of definition of \makeagttitle
\count0=\startpage
\agt\hfill      %   Journal title (top left) 
%   logo (top right)
\hbox to 45truept{\vbox to 0pt{\vglue -13truept{\logomed A\kern -.37em{\logobig 
T}\kern -.38em G}\vss}\hss}
\break
{\small Volume \thevolumenumber\ (\thevolumeyear)
\startpage--\finishpage\nl
Published: \publishdate}

\vglue .25truein

% title
{\parskip=0pt\leftskip 0pt plus
1fil\def\\{\par\smallskip}{\Large\bf\thetitle}\par\medskip} \vglue
0.05truein

% authors :
%
{\parskip=0pt\leftskip 0pt plus 1fil\def\\{\par}{\sc\theauthors}
\par\medskip}%
 
\vglue 0.03truein 

%  abstract and classification numbers:

{\small\leftskip 25truept\rightskip 25truept{\bf Abstract}\stdspace\theabstract

{\bf AMS Classification}\stdspace\theprimaryclass
\ifx\thesecondaryclass\relax\else; \thesecondaryclass\fi\par
{\bf Keywords}\stdspace \thekeywords\par}\vglue 7truept

}   %%%% end of definition of \makeagttitle

\ifplaintex
%  get print centerpage:
\hoffset 14truemm
\voffset 31truemm
%  fonts for headline and footline
\font\phead=cmsl9 scaled 950
\font\pnum=cmbx10 scaled 913
\font\pfoot=cmsl9 scaled 950
%  headline and footline
\headline{\vbox to 0pt{\vskip -4.5mm\line{\small\phead\ifnum
\count0=\startpage ISSN 1472-2739 (on-line) 1472-2747 (printed)
\hfill {\pnum\folio}\else\ifodd\count0\def\\{ }% 
\ifx\theshorttitle\relax\thetitle\else\theshorttitle\fi\hfill{\pnum\folio}
\else\def\\{ and }{\pnum\folio}\hfill\ifx\theshortauthors\relax\theauthors
\else\theshortauthors\fi\fi\fi}\vss}}
\footline{\vbox to 0pt{\vglue 0mm\line{\small\pfoot\ifnum\count0=\startpage
\copyright\ \gtp\hfill\else
\agt, Volume \thevolumenumber\ (\thevolumeyear)\hfill\fi}\vss}}
\else
%  get print centerpage:
\headsep 23pt
\footskip 35pt
\hoffset -4truemm
\voffset 12.5truemm
%  fonts for headline and footline
\font\lhead=cmsl9 scaled 1050
\font\lnum=cmbx10 
\font\lfoot=cmsl9 scaled 1050
\makeatletter
%  headline and footline
\def\@oddhead{{\small\lhead\ifnum\count0=\startpage ISSN 1472-2739 
(on-line) 1472-2747 (printed)\hfill {\lnum\number\count0}\else\ifodd\count0
\def\\{ }\ifx\theshorttitle\relax \thetitle \else\theshorttitle\fi\hfill
{\lnum\number\count0}\else\def\\{ and }{\lnum\number\count0}
\hfill\ifx\theshortauthors\relax 
\theauthors\else\theshortauthors\fi\fi\fi}}\def\@evenhead{\@oddhead}
\def\@oddfoot{\small\lfoot\ifnum\count0=\startpage\copyright\ \gtp\hfill\else
\agt, Volume \thevolumenumber\ (\thevolumeyear)\hfill\fi}
\def\@evenfoot{\@oddfoot}
\makeatother
\fi
%  force \makeagttitle
\let\maketitlepage\makeagttitle
\let\makeshorttitle\maketitlepage
\let\maketitle\maketitlepage

   %%%comment out to create xxx header file

\newwrite\gtoutfile
\long\gdef\makeheadfile{  %%% start of definition of \makeheadfile
{\def\\{, }\def\s{ }
\immediate\openout\gtoutfile head.xxx
\immediate\write\gtoutfile{To: math@arxiv.org}
\immediate\write\gtoutfile{Subject: put OR rep NNNNN:ppppp}
\immediate\write\gtoutfile{--text follows this line--}
\immediate\write\gtoutfile{Proxy-for: \ifx\theasciiauthors\relax
\theauthors\else\theasciiauthors\fi\s<\ifx\theasciiemail\relax\theemail\else\theasciiemail\fi>}
\immediate\write\gtoutfile{\noexpand\\}
\immediate\write\gtoutfile{Authors: \ifx\theasciiauthors\relax
\theauthors\else\theasciiauthors\fi}
{\def\\{ }\immediate\write\gtoutfile{Title: \ifx\theasciititle\relax
\thetitle\else\theasciititle\fi}}
\immediate\write\gtoutfile{Subj-class: GT or SG, GR etc}
\immediate\write\gtoutfile{MSC-class: \theprimaryclass\ifx\thesecondaryclass\relax\else, \thesecondaryclass\fi}
\immediate\write\gtoutfile{Journal-ref: Algebr. Geom. Topol. \thevolumenumber\s
(\thevolumeyear) \startpage-\finishpage}
\immediate\write\gtoutfile{Comments: Published by Algebraic and
Geometric Topology at}
\immediate\write\gtoutfile{\s\s\s  http://www.maths.warwick.ac.uk/agt/AGTVol\thevolumenumber/agt-\thevolumenumber-\thepapernumber.abs.html}
\immediate\write\gtoutfile{\noexpand\\}
\immediate\write\gtoutfile{}
\ifx\theasciiabstract\relax
\immediate\write\gtoutfile{\theabstract}\else
\immediate\write\gtoutfile{\theasciiabstract}\fi
\immediate\write\gtoutfile{}
\immediate\write\gtoutfile{\noexpand\\}
\immediate\write\gtoutfile{}
\immediate\closeout\gtoutfile}}  %%% end of definition of \makeheadfile

\def\maketitlepage{\makeagttitle\makeheadfile}
\let\makeshorttitle\maketitlepage
\let\maketitle\maketitlepage

%%%
%%%  This is agtout.tex.  
%%%
%%%  This the version of  gtoutput.tex  intended to finish formatting
%%%  papers published in Algebriac & Geometric Topology and stored in the
%%%  arXiv.   All versions of  gtoutput.tex  are copyright 
%%%  GT Publications and are to be used _only_ for formatting
%%%  the officially published version of ABT or G&T papers.
%%%
%%%
%%%                                             Colin Rourke  27.102000
%%%
%%%  To create header file  head.xxx  comment out the first \endinput

%  test for latex or plain tex
\def\ifplaintex{\expandafter\ifx\csname documentclass\endcsname\relax}

\def\gtp{{\mathsurround=0pt\it $\cal G\mskip-2mu$eometry \&\ 
$\cal T\!\!$opology $\cal P\!$ublications}}  % GT publications

\def\recd{{\small Received:\qua\receiveddate\ifx\reviseddate\relax
\else\qquad Revised:\qua\reviseddate\fi\par}} 

%  define the various new ingredients of the title page and the data
%  output files

\def\lognumber#1{\def\thelognumber{#1}}
\def\volumenumber#1{\def\thevolumenumber{#1}}
\def\volumeyear#1{\def\thevolumeyear{#1}}
\def\papernumber#1{\def\thepapernumber{#1}}
\def\pagenumbers#1#2{\def\startpage{#1}\def\finishpage{#2}}
\def\published#1{\def\publishdate{#1}}

\def\received#1{\def\receiveddate{#1}}
\def\revised#1{\def\reviseddate{#1}}
\def\accepted#1{\def\accepteddate{#1}}
\def\asciititle#1{\def\theasciititle{#1}}

\def\asciiaddress#1{\def\theasciiaddress{#1}}

\long\def\asciiabstract#1{\long\def\theasciiabstract{#1}}

%  initialise

\let\\\par\let\thelognumber\relax\let\thevolumenumber\relax
\let\thepapernumber\relax\let\thevolumeyear\relax\let\startpage\relax
\let\finishpage\relax\let\publishdate\relax\let\receiveddate\relax
\let\reviseddate\relax\let\accepteddate\relax\let\theasciititle\relax
\let\theasciiauthors\relax\let\theasciiaddress\relax
\let\theasciiabstract\relax

\let\theasciiemail\relax

%%%% fonts for AGT logo:

\ifplaintex
\font\logobig=cmssbx10 scaled 3836
\font\logomed=cmssbx10 scaled 2557
\else
\font\logobig=cmssbx10 scaled 4200
\font\logomed=cmssbx10 scaled 2800
\fi

\long\def\makeagttitle{   %%% start of definition of \makeagttitle
\count0=\startpage
\agt\hfill      %   Journal title (top left) 
%   logo (top right)
\hbox to 45truept{\vbox to 0pt{\vglue -13truept{\logomed A\kern -.37em{\logobig 
T}\kern -.38em G}\vss}\hss}
\break
{\small Volume \thevolumenumber\ (\thevolumeyear)
\startpage--\finishpage\nl
Published: \publishdate}

\vglue .25truein

% title
{\parskip=0pt\leftskip 0pt plus
1fil\def\\{\par\smallskip}{\Large\bf\thetitle}\par\medskip} \vglue
0.05truein

% authors :
%
{\parskip=0pt\leftskip 0pt plus 1fil\def\\{\par}{\sc\theauthors}
\par\medskip}%
 
\vglue 0.03truein 

%  abstract and classification numbers:

{\small\leftskip 25truept\rightskip 25truept{\bf Abstract}\stdspace\theabstract

{\bf AMS Classification}\stdspace\theprimaryclass
\ifx\thesecondaryclass\relax\else; \thesecondaryclass\fi\par
{\bf Keywords}\stdspace \thekeywords\par}\vglue 7truept

}   %%%% end of definition of \makeagttitle

\ifplaintex
%  get print centerpage:
\hoffset 14truemm
\voffset 31truemm
%  fonts for headline and footline
\font\phead=cmsl9 scaled 950
\font\pnum=cmbx10 scaled 913
\font\pfoot=cmsl9 scaled 950
%  headline and footline
\headline{\vbox to 0pt{\vskip -4.5mm\line{\small\phead\ifnum
\count0=\startpage ISSN 1472-2739 (on-line) 1472-2747 (printed)
\hfill {\pnum\folio}\else\ifodd\count0\def\\{ }% 
\ifx\theshorttitle\relax\thetitle\else\theshorttitle\fi\hfill{\pnum\folio}
\else\def\\{ and }{\pnum\folio}\hfill\ifx\theshortauthors\relax\theauthors
\else\theshortauthors\fi\fi\fi}\vss}}
\footline{\vbox to 0pt{\vglue 0mm\line{\small\pfoot\ifnum\count0=\startpage
\copyright\ \gtp\hfill\else
\agt, Volume \thevolumenumber\ (\thevolumeyear)\hfill\fi}\vss}}
\else
%  get print centerpage:
\headsep 23pt
\footskip 35pt
\hoffset -4truemm
\voffset 12.5truemm
%  fonts for headline and footline
\font\lhead=cmsl9 scaled 1050
\font\lnum=cmbx10 
\font\lfoot=cmsl9 scaled 1050
\makeatletter
%  headline and footline
\def\@oddhead{{\small\lhead\ifnum\count0=\startpage ISSN 1472-2739 
(on-line) 1472-2747 (printed)\hfill {\lnum\number\count0}\else\ifodd\count0
\def\\{ }\ifx\theshorttitle\relax \thetitle \else\theshorttitle\fi\hfill
{\lnum\number\count0}\else\def\\{ and }{\lnum\number\count0}
\hfill\ifx\theshortauthors\relax 
\theauthors\else\theshortauthors\fi\fi\fi}}\def\@evenhead{\@oddhead}
\def\@oddfoot{\small\lfoot\ifnum\count0=\startpage\copyright\ \gtp\hfill\else
\agt, Volume \thevolumenumber\ (\thevolumeyear)\hfill\fi}
\def\@evenfoot{\@oddfoot}
\makeatother
\fi
%  force \makeagttitle
\let\maketitlepage\makeagttitle
\let\makeshorttitle\maketitlepage
\let\maketitle\maketitlepage

   %%%comment out to create xxx header file

\newwrite\gtoutfile
\long\gdef\makeheadfile{  %%% start of definition of \makeheadfile
{\def\\{, }\def\s{ }
\immediate\openout\gtoutfile head.xxx
\immediate\write\gtoutfile{To: math@arxiv.org}
\immediate\write\gtoutfile{Subject: put OR rep NNNNN:ppppp}
\immediate\write\gtoutfile{--text follows this line--}
\immediate\write\gtoutfile{Proxy-for: \ifx\theasciiauthors\relax
\theauthors\else\theasciiauthors\fi\s<\ifx\theasciiemail\relax\theemail\else\theasciiemail\fi>}
\immediate\write\gtoutfile{\noexpand\\}
\immediate\write\gtoutfile{Authors: \ifx\theasciiauthors\relax
\theauthors\else\theasciiauthors\fi}
{\def\\{ }\immediate\write\gtoutfile{Title: \ifx\theasciititle\relax
\thetitle\else\theasciititle\fi}}
\immediate\write\gtoutfile{Subj-class: GT or SG, GR etc}
\immediate\write\gtoutfile{MSC-class: \theprimaryclass\ifx\thesecondaryclass\relax\else, \thesecondaryclass\fi}
\immediate\write\gtoutfile{Journal-ref: Algebr. Geom. Topol. \thevolumenumber\s
(\thevolumeyear) \startpage-\finishpage}
\immediate\write\gtoutfile{Comments: Published by Algebraic and
Geometric Topology at}
\immediate\write\gtoutfile{\s\s\s  http://www.maths.warwick.ac.uk/agt/AGTVol\thevolumenumber/agt-\thevolumenumber-\thepapernumber.abs.html}
\immediate\write\gtoutfile{\noexpand\\}
\immediate\write\gtoutfile{}
\ifx\theasciiabstract\relax
\immediate\write\gtoutfile{\theabstract}\else
\immediate\write\gtoutfile{\theasciiabstract}\fi
\immediate\write\gtoutfile{}
\immediate\write\gtoutfile{\noexpand\\}
\immediate\write\gtoutfile{}
\immediate\closeout\gtoutfile}}  %%% end of definition of \makeheadfile

\def\maketitlepage{\makeagttitle\makeheadfile}
\let\makeshorttitle\maketitlepage
\let\maketitle\maketitlepage

\lognumber{24}
\volumenumber{2}
\volumeyear{2002}
\papernumber{24}
\published{25 June 2002}
\pagenumbers{499}{518}
\received{12 March 2002}
\revised {5 June 2002}
\accepted{5 June 2002}

\input rlepsf

\useforwardrefs

\def\rtimes{{\msb \char'157}}

\reflist

\key{Bige} {\bf S. Bigelow}, {\it Braid groups are linear}, J. Amer. Math. Soc.
14 (2001), 471--486. 

\key{CW} {\bf A.\,M. Cohen}, {\bf D.\,B. Wales}, {\it Linearity of Artin groups
of finite type}, Israel J. Math., to appear.

\key{Dig} {\bf F. Digne}, {\it On the linearity of Artin braid groups}, J
Algebra, to appear.  

\key{FN} {\bf E. Fadell}, {\bf L. Neuwirth}, {\it Configuration spaces}, Math.
Scand. 10 (1962), 111--118.

\key{Kra1} {\bf D. Krammer}, {\it The braid group $B_4$ is linear}, Invent.
Math. 142 (2000), 451--486.

\key{Kra2} {\bf D. Krammer}, {\it Braid groups are linear}, Ann. Math. 155
(2002), 131-156.

\key{Law1} {\bf R.\,J. Lawrence}, {\it Homology representations of braid 
groups}, Ph.D. Thesis, University of Oxford, 1989.

\key{Law2} {\bf R.\,J. Lawrence}, {\it Homological representations of the Hecke
algebra}, Commun. Math. Phys. 135 (1990), 141--191.

\key{Par} {\bf L. Paris}, {\it Artin monoids inject in their groups}, Comm.
Math. Helv., to appear.

\key{Sal} {\bf M. Salvetti}, {\it Topology of the complement of real
hyperplanes in ${\Bbb C}^N$}, Invent. Math. 88 (1987), 603--618.

\key{Sal1} {\bf M. Salvetti}, {\it The homotopy type of Artin groups}, Math.
Res. Lett. 1 (1994), 565-577.

\endreflist

\title{A note on the Lawrence--Krammer--Bigelow\\representation}
\asciititle{A note on the Lawrence-Krammer-Bigelow\\representation}

\authors{Luisa Paoluzzi\\Luis Paris}

\address{Laboratoire de Topologie, UMR 5584 du CNRS\\Universit\'e de 
Bourgogne, 9, avenue Alain Savary -- BP 47870\\21078 Dijon CEDEX -- France}
\asciiaddress{Laboratoire de Topologie, UMR 5584 du CNRS\\Universite de 
Bourgogne, 9, avenue Alain Savary - BP 47870\\21078 Dijon CEDEX - France}

\email{paoluzzi@u-bourgogne.fr, lparis@u-bourgogne.fr}

\abstract 
A very popular problem on braid groups has recently been solved by
Bigelow and Krammer, namely, they have found a faithful linear
representation for the braid group $B_n$. In their papers, Bigelow and
Krammer suggested that their representation is the monodromy
representation of a certain fibration. Our goal in this paper is to
understand this monodromy representation using standard tools from the
theory of hyperplane arrangements. In particular, we prove that the
representation of Bigelow and Krammer is a sub-representation of the
monodromy representation which we consider, but that it cannot be the
whole representation.
\endabstract 

\asciiabstract{ 
A very popular problem on braid groups has recently been solved by
Bigelow and Krammer, namely, they have found a faithful linear
representation for the braid group B_n. In their papers, Bigelow and
Krammer suggested that their representation is the monodromy
representation of a certain fibration. Our goal in this paper is to
understand this monodromy representation using standard tools from the
theory of hyperplane arrangements. In particular, we prove that the
representation of Bigelow and Krammer is a sub-representation of the
monodromy representation which we consider, but that it cannot be the
whole representation.}

\primaryclass{20F36}\secondaryclass{52C35, 52C30, 32S22}
\keywords{Braid groups, linear representations, Salvetti complexes}

\makeshorttitle

\section{Introduction}\key{intro}

Consider the ring $R={\Bbb Z}[x^{\pm1},y^{\pm1}]$ of Laurent polynomials in two
variables and the (abstract) free $R$-module:
$$V=\bigoplus_{1\le i<j\le n}R\thinspace e_{ij}$$
For $k\in \{1,\dots,n-1\}$ define the $R$-homomorphism $\rho_k\co V\to V$ by
$$\rho_k(e_{ij})=\cases{
xe_{i-1\,j}+(1-x)e_{ij} & if $k=i-1$ \cr
e_{i+1\,j}-xy(x-1)e_{k\,k+1} & if $k=i<j-1$ \cr
-x^2ye_{k\,k+1} & if $k=i=j-1$ \cr
e_{ij}-y(x-1)^2e_{k\,k+1} & if $i<k<j-1$ \cr
e_{i\,j-1}-xy(x-1)e_{k\,k+1} & if $i<j-1=k$ \cr
xe_{i\,j+1}+(1-x)e_{ij} & if $k=j$ \cr
e_{ij} & otherwise. \cr}$$
The starting point of the present work is the following theorem due to Bigelow
[\ref{Bige}] and Krammer [\ref{Kra1}, \ref{Kra2}].
\proc{Theorem}\key{TK}{\rm (Bigelow, [\ref{Bige}]; Krammer, [\ref{Kra1}, 
\ref{Kra2}])}\qua Let $B_n$ be the braid group on $n$ strings, and let
$\sigma_1,\dots,\sigma_{n-1}$ be the standard generators of $B_n$. Then the 
mapping $\sigma_k\mapsto\rho_k$ induces a well-defined faithful representation 
$\rho\co B_n\to{\rm Aut}_R(V)$. In particular, the braid group $B_n$ is linear.
\endproc
Let $V$ be an $R$-module. A {\sl representation} of $B_n$ on $V$ is a
homomorphism $\rho\co B_n \to{\rm Aut}_R(V)$. By abuse of notation, we may
identify the underlying module $V$ with the representation if no confusion is
possible. Two representations $\rho_1$ and $\rho_2$ on $V_1$ and $V_2$,
respectively, are called {\sl equivalent} if there exist an automorphism
$\nu\co R\to R$ and an isomorphism $f\co V_1 \to V_2$ of abelian groups such
that: 
\items
\itemb $f(\rho_1(b)v)=\rho_2(b)f(v)$ for all $b\in B_n$ and all $v\in V_1$;
\itemb $f(\kappa v)=\nu(\kappa)f(v)$ for all $\kappa\in R$ and all $v\in V_1$.
\enditems  
An {\sl LKB representation} is a representation of $B_n$ equivalent to the one
of Theorem \ref{TK} (LKB stands for Lawrence--Krammer--Bigelow).

Let ${\bf D}$ be a disc embedded in ${\Bbb C}$ such that $1,\dots,n$ lie in the interior of ${\bf D}$ (say ${\bf D}=\{z\in{\Bbb C}\mid |z-(n+1)/2|\le
(n+1)/2\}$), and choose a basepoint $P_b$ on the boundary of ${\bf D}$ (say 
$P_b=(n+1)(1-i)/2$). Define a {\sl fork} to be a tree embedded in ${\bf D}$ 
with four vertices $P_b,p,q,z$ and three edges, and such that
$T\cap\partial{\bf D}=\{P_b\}$, $T\cap\{1,\dots,n\}=\{p,q\}$, and all three 
edges have $z$ as vertex. The LKB representation $V$ defined in [\ref{Kra1}] 
is the quotient of the free $R$-module generated by the isotopy classes of 
forks by certain relations. One can easily verify that these relations are 
invariant by the action of $B_n$, viewed as the mapping class group of ${\bf 
D}\setminus\{1,\dots,n\}$, thus $V$ is naturally endowed with a $B_n$-action. 
Krammer in [\ref{Kra1}] stated that a monodromy representation of $B_n$ on a 
twisted homology, $H_2(F_n;\Gamma_\pi)$, is an LKB representation and referred 
to Lawrence's paper [\ref{Law2}] for the proof. The object of the present 
paper is the study of this monodromy representation on $H_2(F_n;\Gamma_\pi)$. 
Let $R\to{\Bbb C}$ be an embedding. The representation considered by Lawrence 
[\ref{Law1}, \ref{Law2}] is isomorphic to $V\otimes{\Bbb C}$, but her geometric 
construction is slightly different from the construction suggested by Krammer 
[\ref{Kra1}]. In his proof of the linearity of braid groups, Bigelow 
[\ref{Bige}] associated to each fork $T$ an element $S(T)$ of 
$H_2(F_n;\Gamma_\pi)$, and used this correspondence to compute the action of 
$B_n$ on $H_2(F_n;\Gamma_\pi)$. A consequence of his calculation is that 
$H_2(F_n;\Gamma_\pi)\otimes{\Bbb Q}(x,y)$ is isomorphic to $V\otimes{\Bbb 
Q}(x,y)$. 

In [\ref{Dig}] and [\ref{CW}], Digne, Cohen and Wales introduced a new
conceptual approach to the LKB representations based on the theory of root
systems, and extended the results of [\ref{Kra2}] to all spherical type Artin
groups. Using the same approach, linear representations have been defined for 
all Artin groups [\ref{Par}], but it is not known whether the resulting
representations are faithful in the non-spherical case.

The formulae in our definition of the LKB representations are those of
[\ref{Dig}]; the formulae of [\ref{Kra1}], [\ref{Bige}] and [\ref{Law1}] can be
obtained by a change of basis which will be given in Section 5. We choose this
basis because it is the most natural basis in our construction and, as pointed
out before, it has an interpretation in terms of root systems which can be
extended to all Artin groups.

Our goal in this paper is to understand the monodromy action on
$H_2(F_n;\Gamma_\pi)$ using standard tools from the theory of hyperplane
arrangements, essentially the so-called Salvetti complexes. These tools are
especially interesting in the sense that they are less specific to the case
``braid groups'' than the tools of Lawrence, Krammer and Bigelow, and we hope
they will be used in the future for constructing linear representations of
other groups like Artin groups. The main result of the paper is the following:
\proc{Theorem}\key{T12} There is a sub-representation $V$ of
$H_2(F_n;\Gamma_\pi)$ such that:
\items
\item{\rm(i)} $V$ is an LKB representation;
\item{\rm(ii)} $V\neq H_2(F_n;\Gamma_\pi)$ if $n\ge3$;
\item{\rm(iii)} if $V'$ is a sub-representation of $H_2(F_n;\Gamma_\pi)$ and $V'$ 
is an LKB  representation, then $V'\subset V$, if $n\ge4$;
\item{\rm(iv)} $V\otimes{\Bbb Q}(x,y)=H_2(F_n;\Gamma_\pi)\otimes{\Bbb Q}(x,y)$.
\enditems
\endproc
We also prove that $H_2(F_n;\Gamma_\pi)$ is a free $R$-module of rank 
$n(n-1)/2$ and give a basis for $H_2(F_n;\Gamma_\pi)$. Note that (ii) and (iii)
imply that $H_2(F_n;\Gamma_\pi)$ is not an LKB representation if $n\ge 4$. This
fact is still true if $n=3$ but, in this case,
one has two minimal LKB representations in $H_2(F_n;\Gamma_\pi)$. The proof of
this fact is left to the reader. Note also that the equality $V \otimes{\Bbb
Q}(x,y)=H_2(F_n;\Gamma_\pi)\otimes{\Bbb Q}(x,y)$ is already known and can be
found in [\ref{Bige}].

We end this section with a detailed description of the monodromy representation
$H_2(F_n;\Gamma_\pi)$.

For $1\le i<j\le n$, let $H_{ij}$ be the hyperplane of ${\Bbb C}^{\,n}$ with
equation $z_i=z_j$, and let $M'_n={\Bbb C}^{\,n}\setminus(\bigcup_{1\le
i<j\le n}H_{ij})$ denote the complement of these hyperplanes. The symmetric
group $\Sigma_n$ acts freely on $M'_n$ and $B_n$ is the fundamental group of
$M'_n/\Sigma_n=M_n$. By [\ref{FN}], the map $p'\co M'_{n+2}\to M'_n$ which
sends $(z_1,\dots,z_{n+2})$ to $(z_1,\dots,z_n)$ is a locally trivial
fibration. Let
$$L_{1t}=\{ z\in{\Bbb C}^{\,2}\mid z_1=t\},\ L_{2t}=\{ z\in{\Bbb C}^{\,2}\mid
z_2=t\},\quad t=1,\dots,n,$$ 
$$L_3=\{ z\in{\Bbb C}^{\,2}\mid z_1=z_2\}.$$
The fibre of $p'$ at $(1,\dots,n)$ is the complement of the above $2n+1$
complex lines:
$$F'_n={\Bbb C}^{\,2}\setminus(\bigcup_{t=1}^n L_{1t}\cup\bigcup_{t=1}^n
L_{2t}\cup L_3)$$
Let $N_n=M'_{n+2}/(\Sigma_n\times\Sigma_2)$. Then $p'\co M'_{n+2}\to M'_n$
induces a locally trivial fibration $p\co N_n\to M_n$ whose fibre is
$F_n=F'_n/\Sigma_2$.
 
Write $\|z\|=\max\{|z_i|\mid i=1,\dots,n\}$ for $z\in{\Bbb C}^{\,n}$. The map
$s'\co M'_n\to M'_{n+2}$ given by
$$s'(z)=\cases{(z,n+1,n+2) & if $\|z\|\le n$ \cr
\vrule height 12pt depth0pt width 0pt(z,{{n+1}\over {n}}\|z\|,{{n+2}\over {n}}\|z\|) & if $\|z\|\ge n$ \cr}$$
is a well-defined section of $p'$ and moreover induces a section $s\co
M_n\to N_n$ of $p$. So, by the homotopy long exact sequence of $p$, the group
$\pi_1(N_n)$ can be written as a semi-direct product $\pi_1(F_n)\rtimes B_n$.

To construct the monodromy representation we need the following two
propositions whose proofs will be given in Sections \ref{homol} and \ref{act},
respectively.
\proc{Proposition}\key{P13} $H_1(F_n)$ is a free ${\Bbb Z}$-module of rank
$n+1$.                     
\endproc
In fact, we shall see that $H_1(F_n)$ has a natural basis
$\{[a_1],\dots,[a_n],[c_1]\}$. Let $H$ be the free abelian group freely
generated by $\{x,y\}$, let $\pi_0\co H_1(F_n)\to H$ be the homomorphism which
sends $[a_i]$ to $x$ for $i=1,\dots,n$, and $[c_1]$ to $y$, and let
$\pi\co\pi_1(F_n)\to H_1(F_n)\to H$ be the composition of the natural
projection $\pi_1(F_n)\to H_1(F_n)$ with $\pi_0$.
\proc{Proposition}\key{P14}
\items
\item{\rm(i)} The kernel of $\pi$ is invariant for the action of $B_n$. In
particular, the action of $B_n$ on $\pi_1(F_n)$ induces an action of $B_n$ on
$H$.
\item{\rm(ii)} The action of $B_n$ on $H$ is trivial.
\enditems
\endproc   
Let ${\widetilde{F_n}}\to F_n$ be the regular covering space associated to
$\pi$. One has $\pi_1({\widetilde{F_n}})=\ker \pi$, $H$ acts freely and
discontinuously on ${\widetilde{F_n}}$, and ${\widetilde{F_n}}/H=F_n$. The
action of $H$ on ${\widetilde{F_n}}$ endows $H_*({\widetilde{F_n}})$ with a
structure of ${\Bbb Z}[H]$-module. This homology group is called {\sl the
homology of $F_n$ with local coefficients} associated to $\pi$, and is denoted 
by $H_*(F_n;\Gamma_\pi)$.
 
Now, Proposition \ref{P14} implies that the fibration $p\co N_n\to M_n$
induces a representation $\rho_\pi\co\pi_1(M_n)=B_n\to{\rm Aut}_{{\Bbb
Z}[H]}(H_*(F_n;\Gamma_\pi))$, called {\sl monodromy representation} on
$H_*(F_n;\Gamma_\pi)$. In this paper, we shall consider the monodromy
representation $\rho_\pi\co\pi_1(M_n)=B_n\to{\rm Aut}_{{\Bbb
Z}[H]}(H_2(F_n;\Gamma_\pi))$ which is the one referred to by Krammer and
Bigelow.
\section{The Salvetti complex}\key{sal}
An {\sl arrangement of lines} in ${\Bbb R}^2$ is a finite family ${\cal A}$
of affine lines in ${\Bbb R}^2$. The complexification of a line $L$ is
the complex line $L_{\Bbb C}$ in ${\Bbb C}^{\,2}$ with the same equation as
$L$. The {\sl complement of the complexification} of ${\cal A}$ is
$$M({\cal A})={\Bbb C}^{\,2}\setminus(\bigcup_{L\in{\cal A}}L_{\Bbb C}).$$
Let ${\cal A}$ be an arrangement of lines in ${\Bbb R}^2$. Then ${\cal A}$
subdivides ${\Bbb R}^2$ into {\sl facets}. We denote by ${\cal F}({\cal A})$
the set of facets and, for $h=0,1,2$, we denote by ${\cal F}_h({\cal A})$ the
set of facets of dimension $h$. A {\sl vertex} is a facet of dimension $0$, an
{\sl edge} is a facet of dimension $1$, and a {\sl chamber} is a facet of
dimension $2$. We partially order ${\cal F}({\cal A})$ with the relation $F<G$
if $F\subset{\overline{G}}$, where ${\overline{G}}$ denotes the closure of $G$.
 
We now define a CW-complex of dimension $2$, called the {\sl Salvetti
complex} of ${\cal A}$, and denoted by $Sal({\cal A})$. This complex has been
introduced by Salvetti in [\ref{Sal}] in the more general setting of hyperplane
arrangements in ${\Bbb R}^n$, $n$ being any positive integer, and Theorem
\ref{T21}, stated below for the case $n=2$, is proved in [\ref{Sal}] for any 
$n$.
 
To every chamber $C\in{\cal F}_2({\cal A})$ we associate a vertex $w_C$ of
$Sal({\cal A})$. The $0$-skeleton of $Sal({\cal A})$ is $Sal_0({\cal
A})=\{w_C\mid C\in{\cal F}_2({\cal A})\}$.

Let $F\in{\cal F}_1({\cal A})$. There exist exactly two chambers $C,D\in{\cal
F}_2({\cal A})$ satisfying $C,D>F$. We associate to $F$ two oriented $1$-cells
of $Sal({\cal A})$: $a(F,C)$ and $a(F,D)$. The source of $a(F,C)$ is $w_C$ and
its target is $w_D$ while the source of $a(F,D)$ is $w_D$ and its target is
$w_C$ (see Figure \ref{F1}). The $1$-skeleton of $Sal({\cal A})$ is the union
of the $a(F,C)$'s, where $F\in{\cal F}_1({\cal A})$, $C\in{\cal F}_2({\cal A})$
and $F<C$.
\figure
\cl{\relabelbox\small
\epsfxsize 8truecm\epsfbox{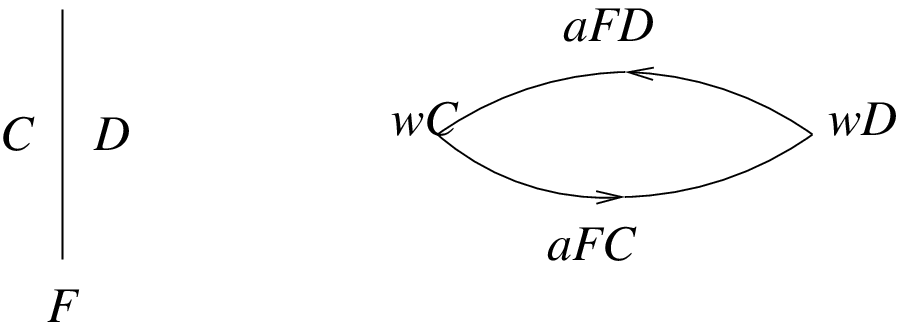}
\adjustrelabel <-2pt,0pt> {C}{$C$}
\relabel {D}{$D$}
\relabel {F}{$F$}
\adjustrelabel <-6pt,-2pt> {wC}{$w_C$}
\adjustrelabel <0pt,-2pt> {wD}{$w_D$}
\relabel {aFC}{$a(F,C)$}
\relabel {aFD}{$a(F,D)$}
\endrelabelbox}
\caption{Edges in $Sal({\cal A})$}
\key{F1}
\endfigure
Let $P\in{\cal F}_0({\cal A})$ and let ${\cal F}_P({\cal A})$ be the set of
chambers $C\in{\cal F}_2({\cal A})$ such that $P<C$. Fix some $C\in{\cal
F}_P({\cal A})$ and write ${\cal F}_P({\cal
A})=\{C,C_1,\dots,C_{n-1},D,D_{n-1},\dots,$ $D_1\}$ (see Figure \ref{F2}). The
set ${\cal F}_P({\cal A})$ has a natural cyclic ordering induced by the
orientation of ${\Bbb R}^2$, so we shall assume the list given above to be
cyclically ordered in this way. Write $C=C_0=D_0$ and $D=C_n=D_n$. For all
$i=1,\dots,n$, there is a unique edge $a_i$ of $Sal_1({\cal A})$ with source
$w_{C_{i-1}}$ and target $w_{C_i}$ and a unique edge $b_i$ of $Sal_1({\cal A})$
with source $w_{D_{i-1}}$ and target $w_{D_i}$. We associate to the pair
$(P,C)$ an oriented $2$-cell $A(P,C)$ of $Sal({\cal A})$ whose boundary is
$$\partial A(P,C)=a_1a_2\dots a_nb_n^{-1}\dots b_2^{-1}b_1^{-1}.$$
The $2$-skeleton of $Sal({\cal A})$ is the union of the $A(P,C)$'s, where
$P\in{\cal F}_0({\cal A})$ and $C\in{\cal F}_P({\cal A})$.
\figure
\cl{\relabelbox\small
\epsfxsize 10truecm\epsfbox{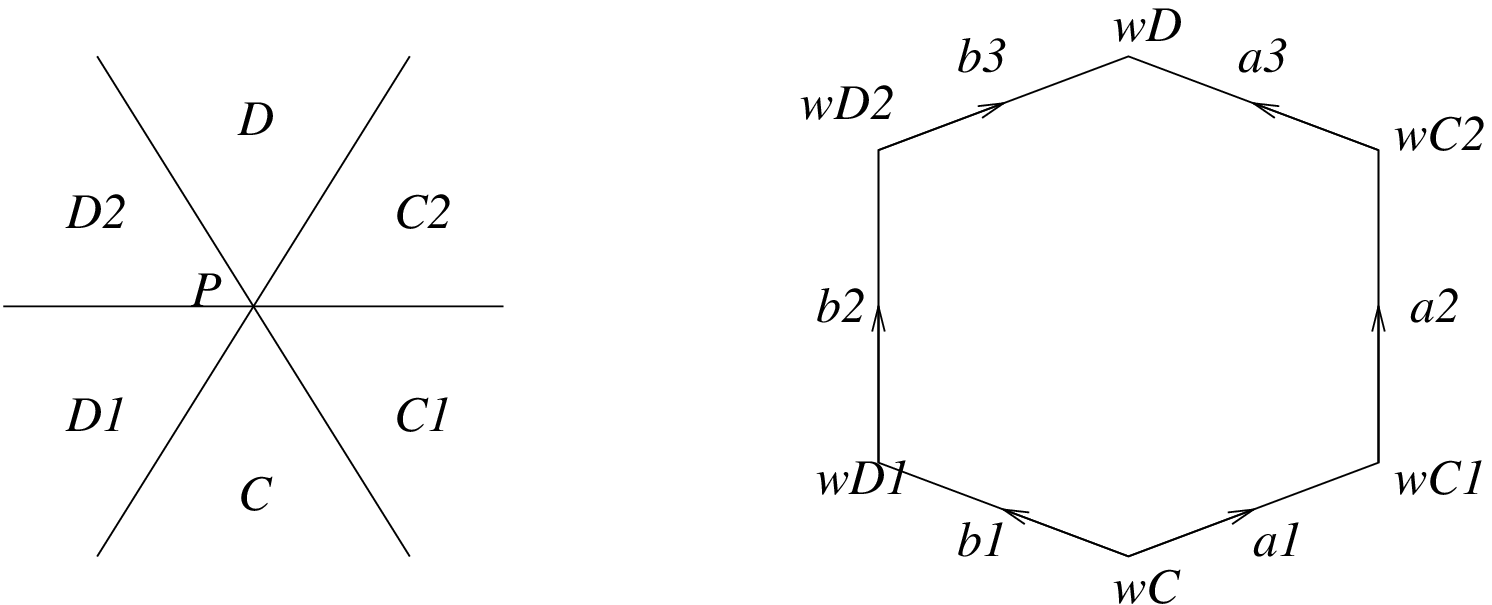}
\relabel {C}{$C$}
\relabel {D}{$D$}
\relabel {C1}{$C_1$}
\relabel {D1}{$D_1$}
\relabel {C2}{$C_2$}
\relabel {D2}{$D_2$}
\adjustrelabel <-4pt, 2pt> {P}{$P$}
\relabel {wC}{$w_C$}
\relabel {wD}{$w_D$}
\relabel {wC1}{$w_{C_1}$}
\relabel {wC2}{$w_{C_2}$}
\relabel {wD1}{$w_{D_1}$}
\relabel {wD2}{$w_{D_2}$} 
\relabel {a1}{$a_1$}
\relabel {a2}{$a_2$}
\relabel {a3}{$a_3$}
\relabel {b1}{$b_1$}
\adjustrelabel <-2pt,0pt> {b2}{$b_2$}
\relabel {b3}{$b_3$}
\endrelabelbox}
\caption{A $2$-cell in $Sal({\cal A})$}
\key{F2}
\endfigure 
\proc{Theorem}\key{T21}{\rm (Salvetti, [\ref{Sal}])}\qua Let ${\cal A}$ be an
arrangement of lines in ${\Bbb R}^2$. There exists an embedding $\delta\co
Sal({\cal A})\to M({\cal A})$ which is a homotopy equivalence.
\endproc
Let ${\cal A}$ be an arrangement of lines in ${\Bbb R}^2$, and let $G$ be a
finite subgroup of ${\rm Aff}({\Bbb R}^2)$ which satisfies:
\items
\itemb $g({\cal A})={\cal A}$ for all $g\in G$;
\itemb $G$ acts freely on ${\Bbb R}^2\setminus(\bigcup_{L\in{\cal A}}L)$.
\enditems 
Then $G$ acts freely on $Sal({\cal A})$ and acts freely on $M({\cal A})$, and
the embedding $\delta\co Sal({\cal A})\to M({\cal A})$ can be chosen to be
equivariant with respect to these actions. Such an equivariant construction can
be found in [\ref{Sal1}] for the particular case where $G$ is a Coxeter group,
and can be carried out in the same way for any group $G$ which satisfies the 
above two conditions. So, $\delta\co Sal({\cal A})\to M({\cal A})$ induces a 
homotopy equivalence ${\bar{\delta}}\co Sal({\cal A})/G\to M({\cal A})/G$.
 
Recall now the spaces $F_n$ and $F'_n$ defined in Section \ref{intro}. Let
$$L_{1\,t}=\{x\in{\Bbb R}^2\mid x_1=t\},\ L_{2\,t}=\{x\in{\Bbb R}^2\mid
x_2=t\},\ t=1,\dots,n,$$
$$L_3=\{x\in{\Bbb R}^2\mid x_1=x_2\},$$
$${\cal A}_n=\{L_{1\,1},\dots,L_{1\,n},L_{2\,1},\dots,L_{2\,n},L_3\}.$$
Then $F'_n=M({\cal A}_n)$ and $F_n=M({\cal A}_n)/\Sigma_2$. The action of 
$\Sigma_2$ on ${\Bbb R}^2$ satisfies:
\items
\itemb $g({\cal A}_n)={\cal A}_n$ for all $g\in\Sigma_2$;
\itemb $\Sigma_2$ acts freely on ${\Bbb R}^2\setminus(\bigcup_{L\in{\cal
A}_n}L)$.
\enditems
It follows that the embedding $\delta\co Sal({\cal A}_n)\to M({\cal A}_n)$
induces a homotopy equivalence ${\bar{\delta}}\co Sal({\cal A}_n)/\Sigma_2\to
M({\cal A}_n)/\Sigma_2=F_n$. 

We now define a new CW-complex, denoted by $Sal(F_n)$, obtained from the
complex $Sal({\cal A}_n)/\Sigma_2$ by collapsing cells, and having the same 
homotopy type as $F_n$. Most of our calculations in Sections \ref{homol} and 
\ref{act} will be based on the description of this complex.

The complex $Sal({\cal A}_n)/\Sigma_2$ can be formally described as follows
(see Figure \ref{F3}):

The set of vertices of $Sal({\cal A}_n)/\Sigma_2$ is 
$$\{P_{ij}\mid1\le i\le j\le n+1\}$$
The set of edges of $Sal({\cal A}_n)/\Sigma_2$ is
$$\{c_i\mid1\le i\le n+1\}\cup\{a_{ij},{\bar{a}}_{ij}\mid1\le i\le j\le
n\}\cup\{b_{ij},{\bar{b}}_{ij}\mid1\le i\le j\le n\}.$$
One has:
$$\vbox{\settabs 2 \columns 
\+ ${\rm source}(a_{ij})={\rm target}({\bar{a}}_{ij})=P_{ij}$ &
${\rm source}(b_{ij})={\rm target}({\bar{b}}_{ij})=P_{i+1\,j+1}$ \cr
\+ ${\rm source}({\bar{a}}_{ij})={\rm target}(a_{ij})=P_{i\,j+1}$ &
${\rm source}({\bar{b}}_{ij})={\rm target}(b_{ij})=P_{i\,j+1}$ \cr
\+ ${\rm source}(c_i)={\rm target}(c_i)=P_{ii}$ & \cr}$$
The set of $2$-cells of $Sal({\cal A}_n)/\Sigma_2$ is
$$\{A_{ijr}\mid1\le i<j\le n\ {\rm and}\ 1\le r\le4\}\cup\{B_{ir}\mid1\le i\le
n\ {\rm and}\ 1\le r\le3\}.$$
One has: 
$$\vbox{\settabs 2 \columns
\+ $\partial A_{ij\,1}=(b_{i\,j-1}a_{ij})(a_{i+1\,j}b_{ij})^{-1}$ &
$\partial B_{i\,1}=
(a_{ii}{\bar{b}}_{ii}c_{i+1})(c_ia_{ii}{\bar{b}}_{ii})^{-1}$ \cr
\+ $\partial A_{ij\,2}=
({\bar{a}}_{i+1\,j}b_{i\,j-1})(b_{ij}{\bar{a}}_{ij})^{-1}$ & 
$\partial B_{i\,2}=
({\bar{b}}_{ii}c_{i+1}b_{ii})({\bar{a}}_{ii}c_ia_{ii})^{-1}$ \cr
\+ $\partial A_{ij\,3}=
(a_{ij}{\bar{b}}_{ij})({\bar{b}}_{i\,j-1}a_{i+1\,j})^{-1}$ &
$\partial B_{i\,3}=
(c_{i+1}b_{ii}{\bar{a}}_{ii})(b_{ii}{\bar{a}}_{ii}c_i)^{-1}$ \cr
\+ $\partial A_{ij\,4}=
({\bar{b}}_{ij}{\bar{a}}_{i+1\,j})({\bar{a}}_{ij}{\bar{b}}_{i\,j-1})^{-1}$ 
\cr}$$
\figure
\cl{\relabelbox\small
\epsfxsize 10truecm\epsfbox{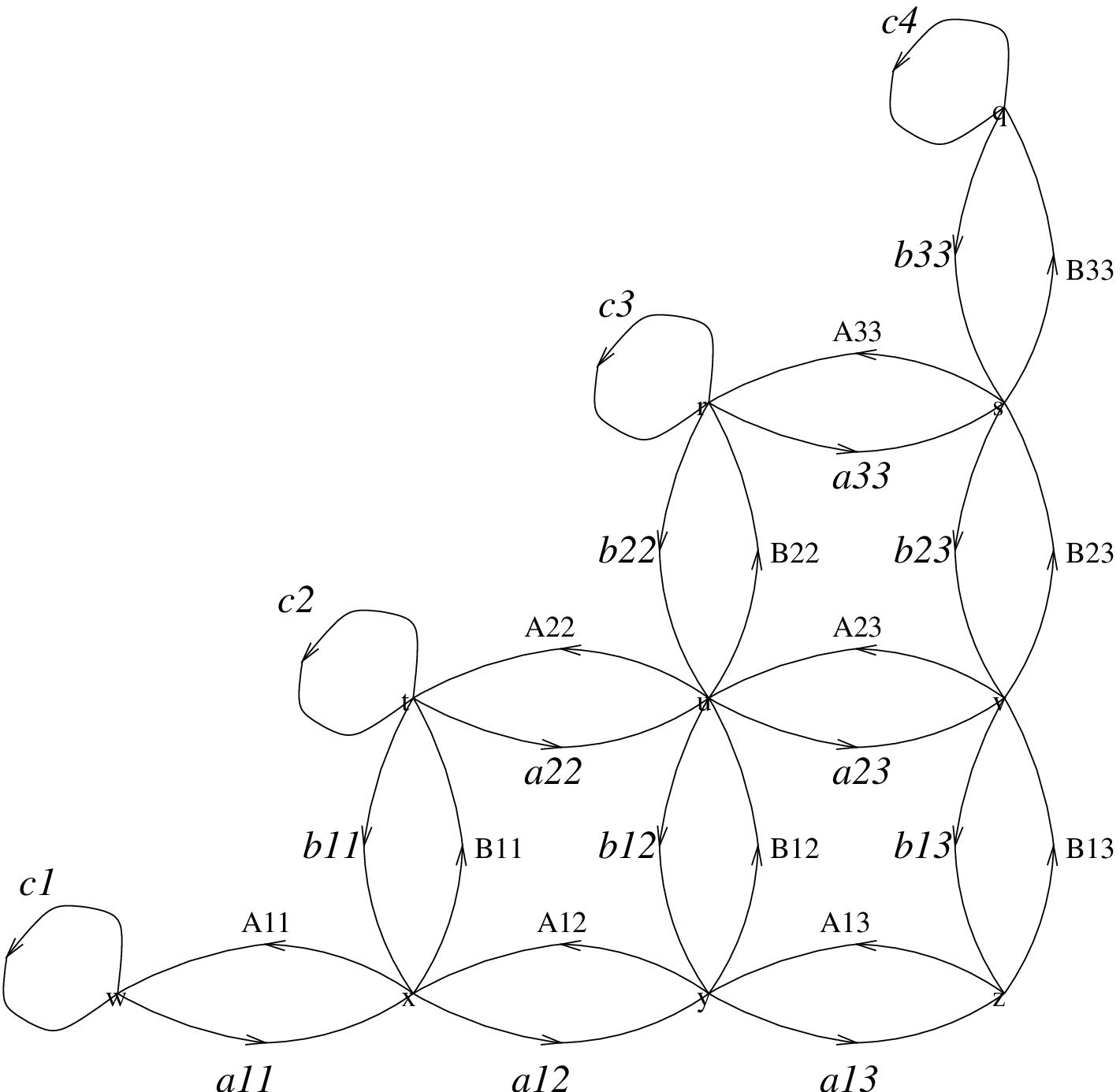}
\adjustrelabel <0pt, 3pt> {a11}{$a_{11}$}
\adjustrelabel <0pt, 3pt> {a12}{$a_{12}$}
\adjustrelabel <0pt, 3pt> {a13}{$a_{13}$}
\relabel {a22}{$a_{22}$}
\relabel {a23}{$a_{23}$}
\relabel {a33}{$a_{33}$}
\relabel {b11}{$b_{11}$}
\relabel {b12}{$b_{12}$}
\relabel {b13}{$b_{13}$}
\relabel {b22}{$b_{22}$}
\relabel {b23}{$b_{23}$}
\relabel {b33}{$b_{33}$}
\relabel {A11}{${\bar{a}}_{11}$}
\relabel {A12}{${\bar{a}}_{12}$}
\relabel {A13}{${\bar{a}}_{13}$}
\relabel {A22}{${\bar{a}}_{22}$}
\relabel {A23}{${\bar{a}}_{23}$}
\relabel {A33}{${\bar{a}}_{33}$}
\relabel {B11}{${\bar{b}}_{11}$}
\relabel {B12}{${\bar{b}}_{12}$}
\relabel {B13}{${\bar{b}}_{13}$}
\relabel {B22}{${\bar{b}}_{22}$}
\relabel {B23}{${\bar{b}}_{23}$}
\relabel {B33}{${\bar{b}}_{33}$}
\relabel {c1}{$c_1$}
\relabel {c2}{$c_2$}
\relabel {c3}{$c_3$}
\relabel {c4}{$c_4$}
\relabel {x}{$\bullet$}
\relabel {y}{$\bullet$}
\relabel {z}{$\bullet$}
\relabel {t}{$\bullet$}
\relabel {w}{$\bullet$}
\relabel {s}{$\bullet$}
\relabel {u}{$\bullet$} 
\relabel {v}{$\bullet$}
\relabel {r}{$\bullet$}
\relabel {q}{$\bullet$}
\endrelabelbox}
\caption{$1$-skeleton of $Sal({\cal A}_3)/\Sigma_2$}
\key{F3}
\endfigure
Let $K$ be the union of all the $A_{ij\,4}$'s. The set $K$ is a subcomplex of
$Sal({\cal A}_n)/\Sigma_2$ which contains all the vertices and all the edges of
$\{{\bar{a}}_{ij},{\bar{b}}_{ij}\mid1\le i\le j\le n\}$, and which is
homeomorphic to a disc. Collapsing $K$ to a single point, we obtain a new
CW-complex denoted by $Sal'(F_n)$. The complex $Sal'(F_n)$ has a unique vertex,
its set of edges is
$$\{c_i\mid1\le i\le n+1\}\cup\{a_{ij},b_{ij}\mid1 \le i\le j\le n\},$$
and its set of 2-cells is
$$\{A_{ijr}\mid1\le i<j \le n\ {\rm and}\ 1\le r\le3\}\cup\{B_{ir}\mid1\le i\le
n\ {\rm and}\ 1\le r\le 3\}.$$
Note that, in $Sal'(F_n)$, the cell $A_{ij\,2}$ is a bigon with boundary
$\partial A_{ij\,2}=b_{i\,j-1}b_{ij}^{-1}$, and $A_{ij\,3}$ is a bigon with
boundary $\partial A_{ij\,3}=a_{ij}a_{i+1\,j}^{-1}$. The complex $Sal(F_n)$ is
obtained from $Sal'(F_n)$ by collapsing all the $A_{ij\,2}$'s for 
$j=i+1,\dots,n$ to a single edge, $b_i=b_{ii}$, and by collapsing all the 
$A_{ij\,3}$'s for $i=1\dots,j-1$ to a single edge, $a_j=a_{jj}$. The complex 
$Sal(F_n)$ has a unique vertex, its set of edges is 
$$\{c_i\mid1\le i\le n+1\}\cup\{a_i,b_i\mid1\le i\le n\},$$
and its set of 2-cells is 
$$\{A_{ij}=A_{ij\,1}\mid 1\le i<j\le n\}\cup\{B_{ir}\mid1\le i\le n\ {\rm
and}\ 1\le r\le3\}.$$
One has: 
$$\vbox{\settabs 2 \columns 
\+ $\partial A_{ij}=(b_ia_j)(a_jb_i)^{-1}$ &
$\partial B_{i\,2}=(c_{i+1}b_i)(c_ia_i)^{-1}$ \cr
\+ $\partial B_{i\,1}=(a_ic_{i+1})(c_ia_i)^{-1}$ &
$\partial B_{i\,3}=(c_{i+1}b_i)(b_ic_i)^{-1}$ \cr}$$
\hfill
\eject
\section{Computing the homology}\key{homol}
For a loop $\alpha$ in $Sal_1(F_n)$, we denote by $[\alpha]$ the element of
$H_1(Sal(F_n))$ represented by $\alpha$. Now, standard methods in homology of
CW-complexes immediately show: 
\proc{Proposition}\key{P31} $H_1(F_n)=H_1(Sal(F_n))$ is the free abelian group
with basis $\{[a_1],\dots,[a_n],[c_1]\}$.
\endproc 
\proc{Remark}\key{R3} We also have the equalities:
$$\eqalign{
[c_i]=[c_1]&\quad{\rm for}\ 1\le i\le n+1\cr
[b_i]=[a_i]&\quad{\rm for}\ 1\le i\le n\cr}$$
\endproc
Recall that $H$ denotes the free abelian group generated by $\{x,y\}$. Define
$\pi_0\co H_1(F_n)\to H$ to be the homomorphism which sends $[a_i]$ to $x$ for
$i=1,\dots,n$, and sends $[c_1]$ to $y$, and let $\pi\co\pi_1(F_n)\to
H_1(F_n)\to H$ be the composition. 

In the following we shall describe a chain complex $C_*(F_n;\Gamma_\pi)$ whose 
homology is $H_*(F_n;\Gamma_\pi)$, define a family $\{ E_{ij}\mid 1\le i<j\le
n\}$ in $H_2(F_n;\Gamma_\pi)$, and prove that $\{ E_{ij}\mid 1\le i<j\le n\}$
is a basis for $H_2(F_n;\Gamma_\pi)\otimes{\Bbb Q}(x,y)$. The sub-module
generated by this family will be the LKB representation $V$ of the statement of
\ref{T12}. We shall end Section \ref{homol} by showing that $H_2(F_n;\Gamma_\pi)$
is a free ${\Bbb Z}[H]$-module.

For $h=0,1,2$, let ${\cal C}_h$ be the set of $h$-cells in $Sal(F_n)$, and let
$C_h(F_n;\Gamma_\pi)$ be the free ${\Bbb Z}[H]$-module with basis ${\cal C}_h$.
Define the differential $d\co C_2(F_n;\Gamma_\pi)\to C_1(F_n;\Gamma_\pi)$
as follows. Let $D\in{\cal C}_2$. Write $\partial
D=\alpha_1^{\varepsilon_1}\dots\alpha_l^{\varepsilon_l}$, where $\alpha_i$ is
an (oriented) 1-cell and $\varepsilon_i \in \{ \pm 1\}$. Set
$$\pi^{(i)}(D)=\cases{\pi(\alpha_1^{\varepsilon_1}\dots
\alpha_{i-1}^{\varepsilon_{i-1}}) & if $\varepsilon_i=1$ \cr
\pi(\alpha_1^{\varepsilon_1}\dots\alpha_{i-1}^{\varepsilon_{i-1}}
\alpha_i^{-1}) & if $\varepsilon_i=-1.$ \cr}$$
Then
$$dD=\sum_{i=1}^l\varepsilon_i\pi^{(i)}(D)\;\alpha_i.$$
The following lemma is a straightforward consequence of this construction.
\proc{Lemma}\key{L32}
\items
\item{\rm(i)} $\ker d=H_2(F_n;\Gamma_\pi)$.
\item{\rm(ii)} Let $d_{\Bbb Q}=d\otimes{\Bbb Q}(x,y)\co
C_2(F_n;\Gamma_\pi)\otimes{\Bbb Q}(x,y)\to C_1(F_n;\Gamma_\pi)\otimes{\Bbb
Q}(x,y)$. Then $\ker d_{\Bbb Q}=H_2(F_n;\Gamma_\pi)\otimes{\Bbb Q}(x,y)$.
\enditems
\endproc
It is easy to obtain the following formulae:
$$\eqalign{
dA_{ij} & =(x-1)(a_j-b_i) \cr
dB_{i\,1} & =(1-y)a_i-c_i+xc_{i+1} \cr
dB_{i\,2} & =-ya_i+yb_i-c_i+c_{i+1} \cr
dB_{i\,3} & =(y-1)b_i-xc_i+c_{i+1} \cr}$$
Now, we define the family $\{E_{ij}\mid1\le i<j\le n\}$. For $1\le i\le n$ 
set:
$$\eqalign{
V_{i\,b} & =-xyB_{i\,1}+x(y-1)B_{i\,2}+B_{i\,3} \cr
V_{i\,a} & =B_{i\,1}+x(y-1)B_{i\,2}-xyB_{i\,3} \cr
V_{i\,0} & =-yB_{i\,1}+(y-1)B_{i\,2}-yB_{i\,3} \cr}$$
For $1\le i<j\le n$ set 
$$E_{ij}=(y-1)(xy+1)A_{ij}+(x-1)V_{i\,b}+(x-1)V_{j\,a}+
\sum_{k=i+1}^{j-1}(x-1)^2V_{k\,0}.$$
The chains $V_{i\,b}$, $V_{i\,a}$ and $V_{i\,0}$ have been found with algebraic
manipulations. Their interest lies in the fact that the support of each of them
is $\{B_{i\,1},B_{i\,2}, B_{i\,3}\}$, the boundary of $V_{i\,b}$ is a multiple
of $c_{i+1}$ minus a multiple of $b_i$, the boundary of $V_{i\,a}$ is a
multiple of $c_i$ minus a multiple of $a_i$, and the boundary of $V_{i\,0}$ is
a multiple of $c_i-c_{i+1}$. More precisely, one has: 
$$\eqalign{dV_{i\,b} & =(y-1)(xy+1)b_i-(x-1)(xy+1)c_{i+1} \cr
dV_{i\,a} & =-(y-1)(xy+1)a_i+(x-1)(xy+1)c_i \cr
dV_{i\,0} & =(xy+1)(c_i-c_{i+1}) \cr}$$
Another fact which will be of importance in our calculations is that all the
$A_{ls}$-coordinates of $E_{ij}$ are zero except the $A_{ij}$-one.
\proc{Proposition}\key{P34}
The set $\{E_{ij}\mid 1\le i<j\le n\}$ is a basis for $\ker d_{\Bbb
Q}=H_2(F_n;\Gamma_\pi)\otimes{\Bbb Q}(x,y)$.
\endproc
\prf 
It is easy to see that $dE_{ij}=0$ for all $1\le i<j\le n$. Moreover, since
$E_{ij}$ is the only element of $\{E_{ls}\mid1\le l<s\le n\}$ such that the
$A_{ij}$-coordinate is nonzero, the set $\{ E_{ij}\mid\le i<j\le n\}$ is
linearly independent.

So, to prove Proposition \ref{P34}, it suffices to show that $\dim(\ker 
d_{\Bbb Q})\le n(n-1)/2$. To do so, we exhibit a linear subspace $W$ of
$C_2(F_n;\Gamma_\pi)\otimes{\Bbb Q}(x,y)$ of codimension $n(n-1)/2$ and prove 
that $d_{\Bbb Q}\vert_W$ is injective.

Let ${\cal B}=\{B_{ir}\mid1\le i\le n\ {\rm and}\ 1\le r\le 3\}$, and let $W$
be the linear subspace of $C_2(F_n;\Gamma_\pi)\otimes{\Bbb Q}(x,y)$ generated
by ${\cal B}$. The codimension of $W$ is clearly $n(n-1)/2$. Let $\eta\co
C_1(F_n;\Gamma_\pi)\otimes{\Bbb Q}(x,y)\to W$ be the linear map defined by 
$$\eqalign{\eta(a_i) & =-(xy-y+1)B_{i\,1}-(y-1)B_{i\,2}+yB_{i\,3} \cr
\eta(b_i) & =-xyB_{i\,1}+x(y-1)B_{i\,2}+B_{i\,3} \cr
\eta(c_i) & =\cases{-y(y-1)B_{i\,1}+(y-1)^2B_{i\,2}-y(y-1)B_{i\,3} & if $1\le
i\le n$ \cr
0 & if $i=n+1$ \cr} \cr}$$
Choose a linear ordering of ${\cal B}$ which satisfies $B_{ir}>B_{i+1\,s}$
for $1\le r,s\le 3$. A straightforward calculation shows that the matrix of
$(\eta \circ d_{\Bbb Q}\vert_W)$ with respect to the ordered basis ${\cal B}$ 
is a triangular matrix with nonzero entries on the diagonal, thus $(\eta\circ 
d_{\Bbb Q}\vert_W)$ is invertible and, therefore, $d_{\Bbb Q}\vert_W$ is
injective.
\endprf 
\proc{Remark}\key{R35}\rm
Let $U \in H_2(F_n;\Gamma_\pi)\otimes{\Bbb Q}(x,y)$. As pointed out before,
$E_{ij}$ is the only element of $\{E_{ls}\mid1\le l<s\le n\}$ such that the
$A_{ij}$-coordinate is nonzero. So, if $\alpha_{ij}$ is the $E_{ij}$-coordinate
of $U$, then $\alpha_{ij}(y-1)(xy+1)$ is the $A_{ij}$-coordinate of $U$.
\endproc
\proc{Proposition}\key{P36}%
$H_2(F_n;\Gamma_\pi)$ is a free ${\Bbb Z}[H]$-module of rank $n(n-1)/2$.
\endproc
\prf
Let 
$$X_{ij}=\cases{E_{ij} & if $j=i+1$ \cr
(E_{1\,2}+E_{2\,3}-E_{1\,3})/(y-1) & if $i=1$ and $j=3$ \cr 
(xyE_{i-1\,i}+(x-1)yE_{i\,i+1} & \cr
\quad -E_{i+1\,i+2}-xyE_{i-1\,i+1} & \cr
\quad +E_{i\,i+2})/(y-1)(xy+1) & if $i\ge 2$ and $j=i+2$ \cr
(E_{i+1\,j-1}-E_{i\,j-1}-E_{i+1\,j} & \cr
\quad +E_{ij})/(y-1)(xy+1) & if $j>i+2$ \cr}$$
and let ${\cal X}=\{X_{ij}\mid1\le i<j\le n\}$. We shall prove that ${\cal X}$
is a ${\Bbb Z}[H]$-basis for $H_2(F_n;\Gamma_\pi)$.

Since $X_{ij}$ is a linear combination (with coefficients in ${\Bbb 
Q}(x,y)$) of $\{E_{ls}\mid1\le l<s\le n\}$, one has $d_{\Bbb Q}X_{ij}=0$. 
Moreover, one can easily verify 
$$X_{ij}=\cases{
(xy+1)A_{1\,2}+(xy+1)A_{2\,3}-(xy+1)A_{1\,3} & \cr
\quad -(x-1)B_{2\,1}+(x^2-1)B_{2\,2}-(x-1)B_{2\,3} & if $i=1$ and $j=3$ \cr
xyA_{i-1\,i}+y(x-1)A_{i\,i+1}-A_{i+1\,i+2} & \cr
\quad -xyA_{i-1\,i+1}+A_{i\,i+2} & \cr 
\quad +x(x-1)B_{i\,2}-(x-1)B_{i\,3} & \cr
\quad -(x-1)B_{i+1\,2}+(x-1)B_{i+1\,3} & if $i\ge 2$ and $j=i+2$ \cr
A_{ij}-A_{i\,j-1}-A_{i+1\,j}+A_{i+1\,j-1} & if $j>i+2$, \cr}$$
thus $X_{ij}\in C_2(F_n;\Gamma_\pi)$. So, $X_{ij}\in H_2(F_n;\Gamma_\pi)$. 

Let $<$ be the linear ordering on $\{A_{ij}\mid1\le i<j\le n\}$ defined by
$A_{ij}<A_{ls}$ if either $j-i<s-l$, or $j-i=s-l$ and $i<l$. The
$A_{ij}$-coordinate of $X_{ij}$ is nonzero and, for $A_{ls}>A_{ij}$, the
$A_{ls}$-coordinate of $X_{ij}$ is zero, thus ${\cal X}$ is linearly
independent. 

It remains to show that any element of $H_2(F_n;\Gamma_\pi)$ can be written as
a linear combination of ${\cal X}$ with coefficients in ${\Bbb Z}[H]$. 
Suppose that there exists $U\in H_2(F_n;\Gamma_\pi)$ which cannot be written as 
a linear combination of ${\cal X}$ with coefficients in ${\Bbb Z}[H]$. 
Write
$$U=\sum\alpha_{ij}A_{ij}+\sum\beta_{ir}B_{ir}$$
where $\alpha_{ij},\beta_{ir}\in{\Bbb Z}[H]$. By Remark \ref{R35}, we have
$$U=\sum{\alpha_{ij}\over(y-1)(xy+1)}E_{ij}.$$
Moreover, $U\neq0$. Let $A_{ij}$ be such that $\alpha_{ij}\neq0$ and
$\alpha_{ls}=0$ for $A_{ls}>A_{ij}$. We choose $U$ so that $A_{ij}$ is minimal
(with respect to the ordering defined above).

Suppose $j>i+2$. The $A_{ij}$-coordinate of $X_{ij}$ is $1$ and, for
$A_{ls}>A_{ij}$, the $A_{ls}$-coordinate of $X_{ij}$ is $0$, thus
$U-\alpha_{ij}X_{ij}$ would contradict the minimality of $A_{ij}$.

Suppose $i\ge2$ and $j=i+2$. Again, the $A_{ij}$-coordinate of $X_{ij}$ is $1$
and, for $A_{ls}>A_{ij}$, the $A_{ls}$-coordinate of $X_{ij}$ is $0$, thus
$U-\alpha_{ij}X_{ij}$ would contradict the minimality of $A_{ij}$.

Suppose $i=1$ and $j=3$. Recall the equality
$U=\sum\alpha_{ls}/((y-1)(xy+1))E_{ls}$. The $B_{n\,1}$-coordinate of $U$ is
$\beta_{n\,1}=(x-1)\alpha_{n-1\,n}/((y-1)(xy+1))$, thus $xy+1$ divides
$\alpha_{n-1\,n}$. For $k=4,\dots,n-1$, the $B_{k\,1}$-coordinate of $U$ is
$\beta_{k\,1}=(x-1)(\alpha_{k-1\,k}-xy\alpha_{k\,k+1})/((y-1)(xy+1))$. It
succesively follows, for $k=n-1,n-2,\dots,4$, that $xy+1$ divides
$\alpha_{k-1\,k}$. The $B_{3,1}$-coordinate, $B_{2,1}$-coordinate, and
$B_{1,1}$-coordinate of $U$ are respectively:
$$\eqalign{\beta_{3\,1} &
=(x-1)(\alpha_{2\,3}+\alpha_{1\,3}-xy\alpha_{3\,4})/((y-1)(xy+1)) \cr
\beta_{2\,1} &
=(x-1)(\alpha_{1\,2}+y\alpha_{1\,3}-xy\alpha_{1\,3}
-xy\alpha_{2\,3})/((y-1)(xy+1)) \cr
\beta_{1\,1} & =-xy(x-1)(\alpha_{1\,2}+\alpha_{1\,3})/((y-1)(xy+1)) \cr}$$
Thus
$$\eqalign{\alpha_{2\,3}+\alpha_{1\,3} & \equiv0\ ({\rm mod}\ xy+1) \cr
\alpha_{1\,2}+\alpha_{2\,3}+(y+1)\alpha_{1\,3} & \equiv0\ ({\rm mod}\ xy+1) \cr
\alpha_{1\,2}+\alpha_{1\,3} & \equiv0\ ({\rm mod}\ xy+1). \cr}$$
Hence $xy+1$ divides $\alpha_{1\,3}$. Let $\alpha_{1\,3}'\in{\Bbb Z}[H]$ be 
such that $\alpha_{1\,3}=(xy+1)\alpha_{1\,3}'$. The $A_{1\,3}$-coordinate of 
$X_{1\,3}$ is $-(xy+1)$ and, for $A_{ls}>A_{1\,3}$, the $A_{ls}$-coordinate of 
$X_{1\,3}$ is $0$, thus $U+\alpha_{1\,3}'X_{1\,3}$ would contradict the 
minimality of $A_{ij}=A_{1\,3}$. 

Suppose $j=i+1$. The $B_{i+1\,1}$-coordinate of $U$ is
$\beta_{i+1\,1}=(x-1)\alpha_{i\,i+1}/((y-1)(xy+1))$, thus $(y-1)(xy+1)$ divides
$\alpha_{i\,i+1}$. Let $\alpha_{i\,i+1}'\in{\Bbb Z}[H]$ such that
$\alpha_{i\,i+1}=(y-1)(xy+1)\alpha_{i\,i+1}'$. The $A_{i\,i+1}$-coordinate
of $X_{i\,i+1}$ is $(y-1)(xy+1)$ and, for $A_{ls}>A_{i\,i+1}$, the
$A_{ls}$-coordinate of $X_{i\,i+1}$ is $0$, thus $U-\alpha_{i\,i+1}'X_{i\,i+1}$
would contradict the minimality of $A_{ij}=A_{i\,i+1}$.
\endprf
\section{Computing the action}\key{act}
We shall see in the next section how to interpret the ``forks'' of Krammer and
Bigelow in our terminology, and, from this interpretation, how to use Bigelow's
calculations [\ref{Bige}, Sec.4] to recover the action of $B_n$ on
$H_2(F_n;\Gamma_\pi)$. In this section, we shall apply our techniques for
calculating the action of $B_n$ on $H_2(F_n;\Gamma_\pi)$. Since most of the
results of the section are well-known, some technical details will be left to
the reader.

Let $k\in\{1,\dots,n-1\}$. Choose some small $\varepsilon>0$ (say
$\varepsilon<1/4$) and an embedding ${\cal V}\co{\bf S}^1\times[0,1]\to{\Bbb 
C}$ which satisfies:
\items
\itemb ${\rm im}{\cal V}=\{z\in{\Bbb C}\; ; 1/2-\varepsilon\le\vert
z-k-1/2\vert\le1/2+\varepsilon\}$;
\itemb ${\cal V}(\zeta,1/2)=k+1/2+\zeta/2$, for all $\zeta\in{\bf S}^1$.
\enditems
Consider the {\sl Dehn twist} $T_k^0\co{\Bbb C}\to{\Bbb C}$ defined by
$$(T_k^0\circ{\cal V})(\zeta,t)={\cal V}(e^{-2i\pi t}\zeta,t)$$
for all $(\zeta,t)\in{\bf S}^1\times[0,1]$, and $T_k^0$ is the identity outside
the image of ${\cal V}$. Note that $T_k^0$ interchanges $k$ and $k+1$ and
fixes the other points of $\{1,\dots,n\}$. Consider now the diagonal
homeomorphism $(T_k^0\times T_k^0)\co{\Bbb C}^{\,2}\to{\Bbb C}^{\,2}$. One has
$(T_k^0\times T_k^0)(F_n')=F_n'$, and $(T_k^0\times T_k^0)$ commutes with the
action of $\Sigma_2$ (namely, $(T_k^0\times T_k^0)\circ g=g\circ(T_k^0\times
T_k^0)$ for all $g\in\Sigma_2$), thus $(T_k^0\times T_k^0)$ induces a
homeomorphism $T_k\co F_n'/\Sigma_2=F_n\to F_n$. Recall that
$\sigma_1,\dots,\sigma_{n-1}$ denote the standard generators of the braid group
$B_n$. Then $T_k$ represents $\sigma_k$, namely
$(T_k)_*=\sigma_k\co\pi_1(F_n)\to\pi_1(F_n)$. 

We assume that $P_b=(n+1,0)$ is the basepoint of $F_n$, we denote by $\ast$ the
unique vertex of $Sal(F_n)$, and choose a homotopy equivalence 
${\bar{\delta}}\co Sal(F_n)\to F_n$ which sends $\ast$ to $P_b$.

Let $\gamma\sim_h\beta$ denote two loops in $F_n$ based at $P_b$ that are
homotopic.
\proc{Lemma}\key{L41} Let $k\in\{1,\dots,n-1\}$. Then
$$\eqalign{T_k({\bar{\delta}}(a_i))\sim_h & \cases{
{\bar{\delta}}(a_{i-1}) & if $k=i-1$ \cr
{\bar{\delta}}(a_ia_{i+1}a_i^{-1}) & if $k=i$ \cr
{\bar{\delta}}(a_i) & otherwise \cr} \cr
T_k({\bar{\delta}}(b_i))\sim_h & \cases{
{\bar{\delta}}(b_ib_{i-1}b_i^{-1}) & if $k=i-1$ \cr
{\bar{\delta}}(b_{i+1}) & if $k=i$ \cr
{\bar{\delta}}(b_i) & otherwise \cr} \cr
T_k({\bar{\delta}}(c_i))\sim_h & \cases{
{\bar{\delta}}(a_{i-1}b_ic_ib_i^{-1}a_{i-1}^{-1}) & if $k=i-1$ \cr
{\bar{\delta}}(c_i) & otherwise. \cr} \cr}$$
\endproc
\prf 
The homotopy relations for ${\bar{\delta}}(a_i)$ follow from the fact that
${\bar{\delta}}(a_i)$ can be drawn in the plane ${\Bbb C}\times\{0\}$ as shown
in Figure \ref{F4}, ${\Bbb C}\times\{0\}$ is invariant by $T_k$, and $T_k$ acts
on ${\Bbb C}\times\{0\}$ as the Dehn twist $T_k^0$. The other homotopy
relations can be proved in the same way.
\endprf
\figure
\relabelbox\small
\epsfxsize 12truecm\epsfbox{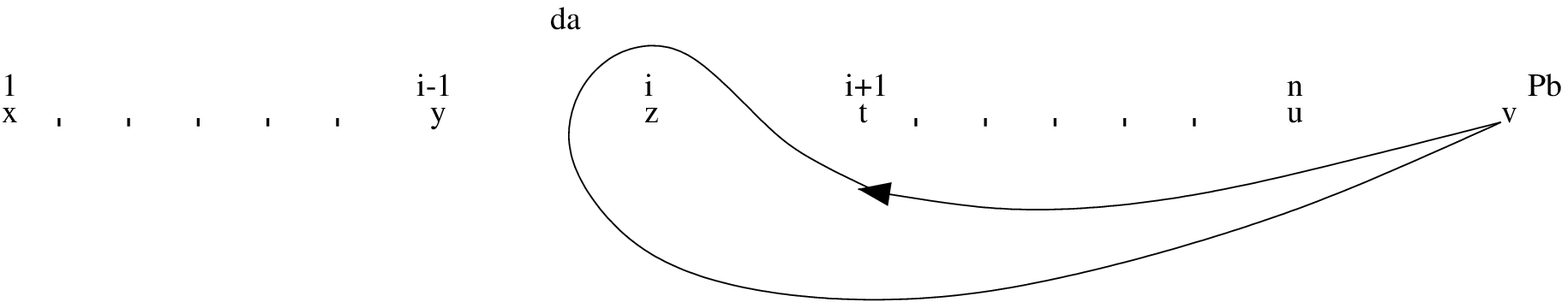}
\relabel {1}{$1$}
\relabel {i-1}{$i-1$}
\relabel {i}{$i$}
\relabel {i+1}{$i+1$}
\relabel {n}{$n$}
\relabel {x}{$\bullet$}
\relabel {y}{$\bullet$}
\relabel {z}{$\bullet$}
\relabel {t}{$\bullet$}
\relabel {u}{$\bullet$}
\relabel {v}{$\bullet$}
\relabel {Pb}{$P_b$}
\relabel {da}{${\bar{\delta}}(a_i)$}
\endrelabelbox
\caption{The curve ${\bar{\delta}}(a_i)$ in ${\Bbb C}\times\{0\}$}
\key{F4}
\endfigure
A straightforward consequence of Lemma \ref{L41} is:
\proc{Corollary}\key{C42} The action of $B_n$ on $H_1(F_n)$ is given by: 
$$\sigma_k([a_i])=\cases{[a_{i-1}] & if $k=i-1$ \cr 
[a_{i+1}] & if $k=i$ \cr
[a_i] & otherwise \cr}$$
$$\sigma_k([c_i])=[c_i]$$
\endproc                          
We now consider the homomorphism $\pi\co\pi_1(F_n)\to H$ defined in Section
\ref{homol}.
\proc{Corollary}\key{C43}
\items
\item{\rm(i)} $\ker\pi$ is invariant by the action of $B_n$. In particular, the
action of $B_n$ on $\pi_1(F_n)$ induces an action of $B_n$ on $H$.
\item{\rm(ii)} The action of $B_n$ on $H$ is trivial.
\enditems
\endproc
So, as pointed out in Section \ref{intro}, this implies:
\proc{Corollary}\key{C44} The locally trivial fibration $p\co N_n\to M_n$
induces a representation $\rho_\pi\co B_n=\pi_1(M_n)\to{\rm Aut}_{{\Bbb
Z}[H]}(H_*(F_n;\Gamma_\pi))$.
\endproc 
We turn now to compute the action of $B_n$ on $H_2(F_n;\Gamma_\pi)$.
 
Let $k\in\{1,\dots,n-1\}$. Define the map $S_k\co Sal_1(F_n)\to Sal(F_n)$ by:
$$\eqalign{
S_k(\ast)= & \ast \cr
S_k(a_i)= & \cases{
a_{i-1} & if $k=i-1$ \cr
a_ia_{i+1}a_i^{-1} & if $k=i$ \cr
a_i & otherwise \cr} \cr
S_k(b_i)= & \cases{
b_ib_{i-1}b_i^{-1} & if $k=i-1$ \cr
b_{i+1} & if $k=i$ \cr
b_i & otherwise \cr} \cr
S_k(c_i)= & \cases{
a_{i-1}b_ic_ib_i^{-1}a_{i-1}^{-1} & if $k=i-1$ \cr
c_i & otherwise \cr} \cr}$$
By Lemma \ref{L41}, $S_k$ induces a homomorphism 
$$(S_k)_*\co\pi_1(Sal(F_n))\to\pi_1(Sal(F_n))$$ 
which is equal to $\sigma_k$. Moreover, by [\ref{FN}],
$Sal(F_n)$ is aspherical, thus $S_k$ extends to a map $S_k\co Sal(F_n)\to
Sal(F_n)$ which is unique up to homotopy.

Let $K$ and $K'$ be two CW-complexes. Call a map $f\co K\to K'$ a {\sl
combinatorial map} if:
\items
\itemb the image of any cell $C$ of $K$ is a cell of $K'$;
\itemb if $\dim C=\dim f(C)$, then $f\vert_C\co C\to f(C)$ is a homeomorphism.
\enditems
We can, and will, suppose that every cell $D$ of $Sal(F_n)$ is endowed with a
cellular decomposition such that $S_k\vert_D\co D\to Sal(F_n)$ is a
combinatorial map. Under this assumption, the map $S_k$ determines a ${\Bbb
Z}[H]$-homomorphism $(S_k)_*\co C_2(F_n;\Gamma_\pi)\to C_2(F_n;\Gamma_\pi)$ as
follows.

Let $D\in{\cal C}_2$ be a $2$-cell of $Sal(F_n)$. Recall that $D$ is endowed
with the cellular decomposition such that $S_k\vert_D\co D\to Sal(F_n)$ is a
combinatorial map. Let ${\cal C}^0_2(D)$ denote the set of $2$-cells $R$ in $D$
such that $S_k(R)$ is a $2$-cell of $Sal(F_n)$. 

Let ${\bf D}^2=\{z\in{\Bbb C}\mid|z|\le1\}$ be the standard disc. In the
definition of the differential $d\co C_2(F_n;\Gamma_\pi)\to
C_1(F_n;\Gamma_\pi)$, for a given 2-cell $D\in{\cal C}_2$, the expression
$\partial D=\alpha_1^{\varepsilon_1}\dots\alpha_l^{\varepsilon_l}$ means that
$D$ is endowed with a cellular map $\phi_D\co{\bf D}^2\to D$ such that the
restriction of $\phi_D$ to the interior of ${\bf D}^2$ is a homeomorphism onto
the interior of $D$, and
$\partial\phi_D=\alpha_1^{\varepsilon_1}\dots\alpha_l^{\varepsilon_l}$, where
$\partial\phi_D\co[0,1]\to D$ is defined by $\partial\phi_D(t)=\phi_D(e^{2i\pi
t})$. Now, every 2-cell $R\in{\cal C}_2^0(D)$ is also endowed with a cellular
map $\phi_R\co{\bf D}^2\to R$ defined by $S_k\circ\phi_R=\phi_{S_k(R)}$. For
$R\in{\cal C}_2^0(D)$, we set $Q_R=\phi_R(1)$. This point should be
understood as the starting point of the reading of the boundary of $R$. 

For every $R\in{\cal C}_2^0(D)$ we set $\epsilon(R)=1$ if $S_k\co R\to
S_k(R)$ preserves the orientation and $\epsilon(R)=-1$ otherwise, and we choose
a path $\gamma_R$ from $\ast$ to $Q_R$ in the 1-skeleton of $D$. Then 
$$(S_k)_*(D)=\sum_{R\in{\cal
C}^0_2(D)}\epsilon(R)(\pi\circ S_k)(\gamma_R)\;S_k(R)$$ 
The sub-module $\ker d=H_2(F_n;\Gamma_\pi)$ is invariant by $(S_k)_*$ and the
restriction of $(S_k)_*$ to $H_2(F_n;\Gamma_\pi)$ is equal to the action of
$\sigma_k$ on $H_2(F_n;\Gamma_\pi)$. 
\proc{Lemma}\key{L45} One can choose $S_k$ such that:  
$$\eqalign{
(S_k)_\ast(A_{ij})= & \cases{
(1-x)A_{ij}+xA_{i-1\,j} & if $k=i-1$ \cr
A_{i+1\,j} & if $k=i<j-1$ \cr
U_i & if $k=i=j-1$ \cr
A_{i\,j-1} & if $i<j-1=k$ \cr
(1-x)A_{ij}+xA_{i\,j+1} & if $k=j$ \cr
A_{ij} & otherwise \cr} \cr
(S_k)_\ast(B_{i\,1})= & \cases{
xB_{i\,3} & if $k=i-1$ \cr
B_{i\,1}+xB_{i+1\,1}-x^2B_{i+1\,3} & if $k=i$ \cr
B_{i\,1} & otherwise \cr} \cr
(S_k)_\ast(B_{i\,2})= & \cases{
U_{i-1}+B_{i\,3}-xB_{i-1\,1}+xB_{i-1\,2} & if $k=i-1$ \cr
B_{i\,1}+xB_{i+1\,2}-xB_{i+1\,3}+yU_i & if $k=i$ \cr
B_{i\,2} & otherwise \cr} \cr
(S_k)_\ast(B_{i\,3})= & \cases{
B_{i\,3}+xB_{i-1\,3}-x^2B_{i-1\,1}-x(y-1)U_{i-1} & if $k=i-1$ \cr
(y-1)U_i+xB_{i\,1} & if $k=i$ \cr
B_{i\,3} & otherwise \cr} \cr}$$
where
$$U_i=(x-1)(B_{i\,1}-B_{i\,2}-B_{i+1\,2}+B_{i+1\,3})-yA_{i\,i+1}.$$
\endproc
\prf 
The method for constructing the extension of $S_k\co Sal_1(F_n)\to Sal(F_n)$ is
as follows. For every $D\in{\cal C}_2$, we compute $S_k(\partial D)$, and, from
this result, we construct a cellular decomposition of $D$ and a combinatorial
map $S_{k,D}\co D\to Sal(F_n)$ which extends the restriction of $S_k$ to
$\partial D$. This can be done case by case without any difficulty. The
maps $S_{k,D}$, $D\in{\cal C}_2$, determine the required extension  $S_k\co
Sal(F_n)\to Sal(F_n)$. With this construction, it is easy to compute
$(S_k)_*(D)$ from the definition given above. 
\endprf
Now, from Lemma \ref{L45} , one can easily compute the action of $B_n$ on
$H_2(F_n;\Gamma_\pi)$ and obtain the following formulae. 
\proc{Proposition}\key{P46} Let $k\in\{1,\dots,n-1\}$ and $1\le i<j\le n$. Then
$$\sigma_k(E_{ij})=\cases{
xE_{i-1\,j}+(1-x)E_{ij} & if $k=i-1$ \cr
E_{i+1\,j}-xy(x-1)E_{k\,k+1} & if $k=i<j-1$ \cr
-x^2yE_{k\,k+1} & if $k=i=j-1$ \cr
E_{ij}-y(x-1)^2E_{k\,k+1} & if $i<k<j-1$ \cr
E_{i\,j-1}-xy(x-1)E_{k\,k+1} & if $i<j-1=k$ \cr
xE_{i\,j+1}+(1-x)E_{ij} & if $k=j$ \cr
E_{ij} & otherwise. \cr}$$
\endproc
Now we can prove our main theorem.
\proof{Proof of Theorem \ref{T12}}
Let $V$ be the ${\Bbb Z}[H]$-submodule generated by $\{E_{ij}\mid 1\le i<j\le
n\}$.
\proof{Proof of (i)}
By Proposition \ref{P34}, $V$ is a free ${\Bbb Z}[H]$-module with basis
$\{E_{ij}\mid
1\le i<j\le n\}$, and, by Proposition \ref{P46}, $V$ is a sub-representation of
$H_2(F_n;\Gamma_\pi)$, and is an LKB representation.
\proof{Proof of (ii)}\kern-1.2pt 
The element $X_{1\,3}$ of the proof of Proposition \ref{P36} lies in 
$H_2(\!F_n;\!\Gamma_\pi)$ but does not lie in $V$.
\proof{Proof of (iii)} 
Let $V'$ be a sub-representation of $H_2(F_n;\Gamma_\pi)$ such that $V'$ is an
LKB representation. By definition, $V'$ is a free ${\Bbb Z}[H]$-module, and
there exist a basis $\{E_{ij}'\mid 1\le i<j\le n\}$ for $V'$ and an
automorphism $\nu\co{\Bbb Z}[H]\to{\Bbb Z}[H]$ such that 
$$\sigma_k(E_{ij}')=\cases{
\nu(x)E_{i-1\,j}'+\nu(1-x)E_{ij}' & if $k=i-1$ \cr
E_{i+1\,j}'-\nu(xy(x-1))E_{k\,k+1}' & if $k=i<j-1$ \cr
-\nu(x^2y)E_{k\,k+1}' & if $k=i=j-1$ \cr
E_{ij}'-\nu(y(x-1)^2)E_{k\,k+1}' & if $i<k<j-1$ \cr
E_{i\,j-1}'-\nu(xy(x-1))E_{k\,k+1}' & if $i<j-1=k$ \cr
\nu(x)E_{i\,j+1}'+\nu(1-x)E_{ij}' & if $k=j$ \cr
E_{ij}' & otherwise. \cr}$$
Note that $V'$ is generated as a ${\Bbb Z}[H]$-module by the $B_n$-orbit of
$E_{1\,2}'$, thus, in order to prove that $V'\subset V$, it suffices to
show that $E_{1\,2}'\in V$. 

For $j=3,\dots,n$, let
$$\eqalign{
F_j= & (xy-1)E_{1\,j}-(xy-1)E_{2\,j}+y(1-x)E_{1\,2} \cr
G_j= & x(x^2y+1)E_{1\,j}+(x^2y+1)E_{2\,j}+x^2y(1-x)E_{1\,2}. \cr}$$
It is easy to see that: 
$$\vbox{\settabs 5 \columns
\+ $\sigma_1(E_{1\,2})=-x^2yE_{1\,2}$ & & & & \cr
\+ $\sigma_1(F_j)=-xF_j$ & & for $3\le j\le n$ & & \cr
\+ $\sigma_1(G_j)=G_j$ & & for $3\le j\le n$ & & \cr 
\+ $\sigma_1(E_{ij})=E_{ij}$ & & for $3\le i<j\le n$ & & \cr}$$
The set $\{E_{1\,2}\}\cup\{F_j, G_j\mid 3\le j\le n\}\cup\{E_{i\,j}\mid 3\le
i<j\le n\}$ is a basis for $H_2(F_n;\Gamma_\pi)\otimes{\Bbb Q}(x,y)$, thus the
eigenvalues of $\sigma_1$ are $-x^2y$ of multiplicity $1$, $-x$ of multiplicity
$(n-2)$, and $1$ of multiplicity $(n-1)(n-2)/2$. The same argument shows that
$-\nu(x^2y)$ is an eigenvalue of $\sigma_1$ of multiplicity $1$ and $E_{1\,2}'$
is an eigenvector associated to this eigenvalue. Since $n\ge4$, it follows that
$\nu(-x^2y)=-x^2y$ and $E_{1\,2}'$ is a multiple of $E_{1\,2}$.

Write $E_{1\,2}'=\lambda E_{1\,2}$, where $\lambda\in{\Bbb Q}(x,y)$. The
$A_{1\,2}$-coordinate of $E_{1\,2}'$, and $B_{1\,3}$-coordinate of $E_{1\,2}'$
are $\lambda(y-1)(xy+1)$ and $\lambda(x-1)$, respectively, and both coordinates
lie in ${\Bbb Z}[H]$, thus $\lambda\in{\Bbb Z}[H]$ and $E_{1\,2}'=\lambda 
E_{1\,2}\in V$.
\proof{Proof of (iv)}
This follows directly from Proposition \ref{P34}.
\endprf
\section{Computing the action with the forks}\key{for}
Recall that ${\widetilde{F_n}}\to F_n$ denotes the regular covering associated 
to $\pi\co\pi_1(F_n)\to H$. Choose some $\varepsilon>0$ (say 
$\varepsilon<1/4$), and, for $p=1,\dots,n$, write $v(p)=\{z\in{\Bbb 
C}\mid|z-p|<\varepsilon\}$. Let $T$ be a fork with vertices $P_b,p,q,z$. Let 
$U(T)$ be the set of pairs $\{x,y\}$ in $F_n$ such that either $x$ or $y$ lies 
in $v(p)\cup v(q)$, and let ${\tilde{U}}(T)$ be the pre-image of $U(T)$ in 
$\widetilde{F_n}$. Bigelow [\ref{Bige}] associated to any fork $T$ a disc 
$S^{(1)}(T)$ embedded in $\widetilde{F_n}$ whose boundary lies in 
${\tilde{U}}(T)$, and proved that $(x-1)^2(xy+1)\,\partial S^{(1)}(T)$ is a
boundary in ${\tilde{U}}(T)$. Thus, there exists an immersed surface 
$S^{(2)}(T)$ whose boundary is equal to $(x-1)^2(xy+1)\,\partial S^{(1)}(T)$. 
Note that the surface $S^{(2)}(T)$ is unique since $H_2({\tilde{U}}(T),{\Bbb 
Z})=0$. So, the element $S(T)=S^{(2)}(T)-(x-1)^2(xy+1)\,S^{(1)}(T)$ is a 
well-defined 2-cycle which represents a non-trivial element of 
$H_2({\widetilde{F_n}};{\Bbb Z})=H_2(F_n;\Gamma_\pi)$. Moreover, the mapping 
$T\mapsto S(T)$ is equivariant by the action of $B_n$. 

In [\ref{Kra1}], Krammer defined a family ${\cal T}=\{T_{pq}\mid1\le p<q\le
n\}$ of forks, called {\sl standard forks}, proved that ${\cal T}$ is a basis
for the LKB representation $V$, defined as a quotient of the free ${\Bbb
Z}[H]$-module generated by the isotopy classes of forks, and explicitly
computed the action of $B_n$ on ${\cal T}$. Let ${\widetilde{Sal}}(F_n)\to
Sal(F_n)$ be the regular covering associated to $\pi\co\pi_1(Sal(F_n))\to H$.
Let $T_{pq}$ be a standard fork with vertices $P_b,p,q,z$, and assume $p+1<q$.
Let
$$Sal(T_{pq})=(\cup_{r=1}^3B_{pr})\cup(\cup_{r=1}^3B_{qr})\cup
A_{pq}\cup(\cup_{k=p+1}^{q-1}A_{kq})\cup(\cup_{k=p+1}^{q-1}A_{pk}),$$
and let ${\widetilde{Sal}}(T_{pq})$ be the pre-image of $Sal(T_{pq})$ in
${\widetilde{Sal}}(F_n)$. Then the pair
$({\widetilde{Sal}}(F_n),{\widetilde{Sal}}(T_{pq}))$ is homotopy equivalent to
$({\widetilde{F_n}},{\tilde{U}}(T_{pq}))$. Let 
$$X_{pq}^{(1)}=\cup_{k=p+1}^{q-1}x^{k-1}B_{k\,1}.$$
The set $X_{pq}^{(1)}$ is a disc embedded in ${\widetilde{Sal}}(F_n)$ whose
boundary lies in ${\widetilde{Sal}}(T_{pq})$, and one can choose the homotopy
equivalence 
$$({\widetilde{Sal}}(F_n),{\widetilde{Sal}}(T_{pq}))\to
({\widetilde{F_n}},{\tilde{U}}(T_{pq}))$$ 
such that $X_{pq}^{(1)}$ is sent to $S^{(1)}(T_{pq})$. Let 
$$X_{pq}^{(2)}=
x^p(x-1)V_{p\,b}+x^{q-1}(x-1)V_{q\,a}+x^{q-1}(y-1)(xy+1)A_{pq}+$$
$$-\sum_{k=p+1}^{q-1}x^{k-1}(x-1)(y-1)(xy+1)A_{pk}.$$
$X_{pq}^{(2)}$ is a 2-chain in $C_2({\widetilde{Sal}}(T_{pq});{\Bbb Z})$ and 
one has 
$$dX_{pq}^{(2)}=(x-1)^2(xy+1)\,dX_{pq}^{(1)}=$$
$$=(x-1)^2(xy+1)\left(-\sum_{k=p+1}^{q-1}x^{k-1}(y-1)a_k 
-x^p c_{p+1}+x^{q-1}c_q\right).$$
(Here, $X_{pq}^{(1)}=\sum_{k=p+1}^{q-1}x^{k-1}B_{k\,1}$ is viewed as a
2-chain). In particular, one has the equality
$S(T_{pq})=X_{pq}^{(2)}-(x-1)^2(xy+1)X_{pq}^{(1)}$ in
$H_2(F_n;\Gamma_\pi)=H_2(Sal(F_n);\Gamma_\pi)$. 

A similar argument shows that $S(T_{pq})=x^pE_{pq}$ if $q=p+1$. 

Let 
$$X_{pq}=\cases{
x^pE_{pq} & if $q=p+1$ \cr
X_{pq}^{(2)}-(x-1)^2(xy+1)X_{pq}^{(1)} & if $q>p+1$. \cr}$$
Note that, by Remark \ref{R35}, one has
$$X_{pq}=x^{q-1}E_{pq}-\sum_{k=p+1}^{q-1} x^{k-1}(x-1)E_{pk}.$$
The set $\{X_{pq}\mid1 \le p<q\le n\}$ is a ${\Bbb Z}[H]$-basis for the LKB
representation $V$ of the statement of \ref{T12} and, by the above remarks, the
action of $B_n$ on the $X_{pq}$'s is given by the formulae of [\ref{Bige}, Thm.
4.1].

\references\Addresses\recd
 
\bye